\documentclass[a4paper,11pt]{article}
\usepackage{latexsym,amsfonts,amsmath,amssymb}
\usepackage{graphics, graphicx}
\usepackage{vmargin}

\numberwithin{equation}{section}
\numberwithin{figure}{section}
\numberwithin{table}{section}

\newcommand {\nl}{\newline}

\newcommand {\ph}[2]{\phi(#1,#2)}

\newcommand {\tstar}[1]{t^{*}(#1)}
\newcommand {\beq}{\begin{eqnarray*}}
\newcommand {\eeq}{\end{eqnarray*}}
\newcommand {\e}{\mathrm{e}}
\newcommand {\expl}[2]{\mathrm{e}^{-\Lambda(#1,#2)}}
\newcommand {\expla}[2]{\mathrm{e}^{-\alpha #2-\Lambda(#1,#2)}}
\newcommand{\h}[1]{\widehat{#1}}
\newcommand {\Li}[1]{\big[ #1 \big]}
\newcommand{\widebar}[1]{\overline{#1}}

\def\ds{\displaystyle}
\def\1{\mathbf{1}}

\def\NN{\mathbb{N}}

\def\RR{\mathbb{R}}

\begin{document}
\newtheorem{theorem}{Theorem}[section]
\newtheorem{proposition}[theorem]{Proposition}
\newtheorem{lemma}[theorem]{Lemma}
\newtheorem{corollary}[theorem]{Corollary}
\newtheorem{definition}[theorem]{Definition}
\newtheorem{remark}[theorem]{Remark}
\newtheorem{conjecture}[theorem]{Conjecture}
\newtheorem{assumption}[theorem]{Assumption}

\title{Numerical method for impulse control of Piecewise Deterministic Markov Processes\footnote{This work was supported by ARPEGE program of the French National Agency of Research (ANR),
project ''FAUTOCOES'', number ANR-09-SEGI-004.}}
\author{B. de Saporta and F. Dufour\\Universit\'e de Bordeaux, IMB, CNRS, UMR 5251,\\ and INRIA Bordeaux, team CQFD,\\
351 cours de la Lib\'eration\\
F-33405 Talence Cedex, France\\
\small saporta@math.u-bordeaux1.fr, dufour@math.u-bordeaux1.fr}
\date{\today}

\maketitle

\begin{abstract}
This paper presents a numerical method to calculate the value function for a general discounted  impulse control problem for piecewise deterministic Markov processes.
Our approach is based on a quantization technique for the underlying Markov chain defined by the post jump location and inter-arrival time. Convergence results are obtained and more importantly we are able to give a convergence rate of the algorithm.
The paper is illustrated by a numerical example.

\medskip

\noindent\textbf{Keywords}: Impulse control, Piecewise deterministic Markov processes, Long-run average cost, Numerical approximation, Quantization.

\medskip

\noindent\textbf{MSC classification}: 93E25, 60J25, 93E20, 93C57.
\end{abstract}


\section{Introduction}
\label{SecIntro}
We present here a numerical method to compute the value function of an impulse control problem for a piecewise deterministic Markov process. Our approach is based on the quantization of an underlying discrete-time Markov chain related to the continuous-time process and path-adapted time discretization grids.

\bigskip

Piecewise-deterministic Markov processes (PDMP's) have been introduced in the literature by M. Davis \cite{davis93} as a general class of stochastic hybrid models.
PDMP's are a family of Markov processes involving deterministic motion punctuated by random jumps.
The motion of the PDMP includes both continuous and discrete variables $\{(X(t),\Upsilon(t))\}$. 
The hybrid state space (continuous/discrete) is defined as $\RR^{d}\times M$ where $M$ is a countable set. 
The process depends on three local characteristics, namely the flow $\phi$, the jump rate $\lambda$ and the transition measure $Q$, which specifies the post-jump location.
Starting from $(x,\nu)\in \RR^{d}\times M$ the motion of the process follows the trajectory $(\phi_{\nu}(x,t),\nu)$ until the first jump time $T_1$ which occurs either spontaneously in a Poisson-like fashion with rate $\lambda_{\nu}(\phi_{\nu}(x,t))$ or when the flow $\phi_{\nu}(x,t)$ hits the boundary of the state-space.
In either case the location of the process at the jump time $T_1$: $\big(X(T_{1}),\Upsilon(T_{1})\big)=\big(Z_{1}, y_{1}\big)$ is selected by the transition measure $Q_{\nu}(\phi_{\nu}(x,T_1),\cdot)$.
Starting from $\big(Z_{1}, y_{1}\big)$, we now select the next inter-jump time $T_{2}-T_{1}$ and postjump location $\big(X(T_{2}),\Upsilon(T_{2})\big)=\big(Z_{2},y_{2}\big)$.
This gives a piecewise deterministic trajectory for $\{(X(t),\Upsilon(t))\}$ with jump times $\{T_{k}\}$ and post jump locations $\{(Z_{k},y_{k})\}$ which follows the flow $\phi$ between two jumps.
A suitable choice of the state space and the local characteristics $\phi$, $\lambda$, and $Q$ provides stochastic models covering a great number of problems of operations research, see \cite{davis93}.
To simplify notation, there is no loss of generality in considering that the state space of the PDMP is taken simply as a subset of
$\RR^{d}$ rather than a product space $\RR^{d}\times M$
as described above, see Remark 24.9 in \cite{davis93} for details.

\bigskip

An impulse control strategy consists in a sequence of single interventions introducing a jump of the process at some controller-specified stopping time and moving the process at that time to some new point in the state space. Our impulse control problem consists in choosing a strategy (if it exists) that minimizes the expected sum  of discounted running and intervention costs up to infinity, and computing the optimal cost thus achieved. Many applied problems fall into this class, such as inventory problems in which a sequence of restocking decisions is made, or optimal maintenance of complex systems with components subject to failure and repair.

\bigskip

Impulse control problems of PDMP's in the context of an expected discounted cost have been considered in
\cite{CD89,dempster_ye95,gatarek91,gatarek92,lenhart89}.
Roughly speaking, in \cite{CD89} the authors study this impulse control problem by using the value improvement approach while in
\cite{dempster_ye95,gatarek91,gatarek92,lenhart89} the authors choose to analyze it by using the variational inequality approach.
In \cite{CD89}, the authors also consider a numerical procedure.
By showing that iteration of the single-jump-or-intervention operator generates a sequence of functions converging to the value function of the problem, they derive an algorithm to compute an approximation of that value function. 
Their approach is also based on a uniform discretization of the state space similar to the one proposed by H. J. Kushner in \cite{kushner77}. 
In particular, they derive a convergence result for the approximation scheme but no estimation of the rate of convergence is given.
To the best of our knowledge, it is the only paper presenting a computational method for solving the impulse control problem for a PDMP in the context of discounted cost. Remark that a similar procedure has been applied by O. Costa in \cite{costa93} to derive a numerical scheme for the impulse control problem with a long run average cost.

\bigskip 

Our approach is also based on the iteration of the single-jump-or-intervention operator, but we want to derive a convergence rate for our approximation. 
Our method does not rely on a \emph{blind} discretization of the state space, but on a discretization that depends on time and takes into account the random nature of the process.
Our approach involves a quantization procedure. 
Roughly speaking, quantization is a technique that approximates a continuous state space random variable $X$ by a a random variable $\widehat{X}$ taking only finitely many values and such that the  difference between $X$ and $\widehat{X}$ is minimal for the $L_p$ norm.
Quantization methods have been developed recently in numerical probability, nonlinear filtering or optimal stochastic control with applications in finance, see e.g. \cite{bally03,bally05,pages98,pages05,pages04b,pages04} and references therein.
It has also been successfully used by the authors to compute an approximation of the value function and optimal strategy for the optimal stopping problem for PDMP's \cite{AAP10}.

\bigskip 

Although the value function of the impulse control problem can be computed by iterating implicit optimal stopping problems, see \cite{CD89} Proposition~2 or \cite{davis93} Proposition 54.18, from a numerical point of view the impulse control is much more difficult to handle than the optimal stopping problem.
Indeed, for the optimal stopping problem, the value function is computed as the limit of a sequence $(v_n)$ constructed by iterating an operator $L$. This iteration procedure yields an iterative construction of a sequence of random variables $v_n(Z_n)$ (where $(Z_n)$  is an embedded discrete-time process).
This was the keystone of our approximation procedure.
As regards impulse control, the iterative construction for the corresponding random variables does not hold anymore, see Section~\ref{section impulse} for details.
This is mostly due to the fact that not only does the controller choose times to stop the process, but they also choose a new starting point for the process to restart from after each intervention.
This makes the single-jump-or-intervention operator significantly more complicated to iterate that the single-jump-or-stop operator used for optimal stopping.
We manage to overcome this extra difficulty by using two series of quantization grids instead of just the one we used for optimal stopping.

\bigskip

The paper is organized as follows. In Section~\ref{SecDef} we give a precise definition of a PDMP and state our notation and assumptions. In Section~\ref{section impulse}, we present the impulse control problem and recall the iterative construction of the value function presented in \cite{CD89}. In Section~\ref{ApproxValue}, we explain our approximation procedure and prove its convergence with error bounds.
Finally in Section~\ref{section appli} we present a numerical example. Some technical results are postponed to the Appendix.

\section{Definitions and assumptions}
\label{SecDef}
We first give a precise definition of a piecewise deterministic Markov process (PDMP).
Some general assumptions are presented in the end of this section.
Let us introduce first some standard notation. Let $M$ be a metric space. $\mathbf{B}(M)$ is the set of real-valued, bounded, measurable functions defined on $M$. The Borel $\sigma$-field of $M$ is denoted by
$\mathcal{B}(M)$. Let $Q$ be a Markov kernel on $(M,\mathcal{B}(M))$ and $w\in \mathbf{B}(M)$, $\ds Qw(x)=\int_{M} w(y) Q(x,dy)$ for $x\in M$.
For $(a,b)\in \RR^2$, $a\wedge b= \min(a,b)$ and $a\vee b= \max(a,b)$.

Let $E$ be an open subset of $\RR^d$, $\partial E$ its boundary and $\widebar{E}$ its closure. 
A PDMP is determined by its local characteristics $(\phi,\lambda,Q)$ where:
\nl
$\bullet$ the flow $\phi: \RR^{d} \times \RR \to \RR^{d}$ is a one-parameter group of homeomorphisms: $\phi$ is continuous,
$\phi(\cdot , t)$ is an homeomorphism for each $t\in \RR$ satisfying $\phi(\cdot , t+s)=\phi(\phi(\cdot , s),t)$.
\nl
For all $x$ in $E$, let us denote
$$t^{*}(x)\doteq \inf \{t>0:\phi(x,t)\in \partial E \},$$
with the convention $\inf \emptyset = \infty$.
\nl
$\bullet$ the jump rate $\lambda: \widebar{E} \to \RR_{+}$ is assumed to be a measurable function.
\nl
$\bullet$ $Q$ is a Markov kernel on $(\widebar{E},\mathcal{B}(\widebar{E}))$ satisfying the following property:
$$(\forall x\in \widebar{E}), \quad Q(x,E-\{x\})=1.$$

\noindent
From these characteristics, it can be shown \cite[p. 62-66]{davis93} that there exists a filtered probability space $(\Omega,\mathcal{F},\{ \mathcal{F}_{t} \}, \{ \mathbf{P}_{x} \}_{x\in E})$
such that the motion of the process $\{X(t)\}$ starting from a point $x\in E$ may be constructed as follows.
Take a random variable $T_1$ such that
\begin{equation*}
\mathbf{P}_x(T_1>t) \doteq
\begin{cases}
e^{-\Lambda(x,t)} & \text{for } t<t^{*}(x),\\
0 & \text{for } t\geq t^{*}(x),
\end{cases}
\end{equation*}
where for $x\in E$ and $t\in [0,t^{*}(x)]$
\begin{equation*}
\Lambda(x,t) \doteq \int_0^t\lambda(\phi(x,s))ds.
\end{equation*}
If $T_1$ generated according to the above probability is equal to infinity, then for $t\in \RR_+$, $X(t)=\phi(x,t)$.
Otherwise select independently an $E$-valued random variable (labelled $Z_{1}$) having distribution $Q(\phi(x,T_1),\cdot )$,
namely $\mathbf{P}_x(Z_{1}\in A)=Q(\phi(x,T_1),A)$ for any $A\in \mathcal{B}(\widebar{E})$.
The trajectory of $\{X(t)\}$ starting at $x$, for $t\leq T_1$ , is given by
\begin{equation*}
X(t) \doteq
\begin{cases}
\phi(x,t)&\text{for } t<T_1, \\
Z_1&\text{for }t=T_1.
\end{cases}
\end{equation*}
Starting from $X(T_1)=Z_1$, we now select the next inter-jump time $T_2-T_1$
and post-jump location $X(T_2)=Z_2$ is a similar way.

\bigskip

This gives a strong Markov process $\{X(t)\}$ with jump times $\bigl\{T_k\bigr\}_{k \in \NN}$ (where $T_{0}= 0$).
Associated to $\{X(t)\}$, there exists a discrete time process $\big( \Theta_{n} \big)_{n\in \NN}$ defined by
$\Theta_{n}=(Z_{n},S_{n})$ with $Z_{n}=X(T_{n})$ and $S_{n}=T_{n}-T_{n-1}$ for $n\geq 1$ and $S_{0}=0$.
Clearly, the process $(\Theta_{n})_{n\in\NN}$ is a Markov chain, and it is the only source of randomness of the process.

\bigskip

We define the following space of functions continuous along the flow with limit towards the boundary:
\begin{align*}
\mathbf{C} =  \bigl\{ & w\in\mathbf{B}(E) \: : \: w(\phi(x,\cdot)): [0,t^{*}(x)) \mapsto \RR \text{ is continuous for each } x\in E\\
& \text{ and whenever } t^{*}(x)< \infty \text{ the limit }\lim_{t\rightarrow t^{*}(x)} w(\phi(x,t)) \text{ exists} \bigr\}.
\end{align*}
For $w\in \mathbf{C}$, we define $\ds w(\phi(x,t^{*}(x)))$ by the limit $\ds \lim_{t\rightarrow t^{*}(x)} w(\phi(x,t))$
(note that the limit exists by assumption). 
Let us introduce $\mathbf{L}$ as the set of functions $w\in \mathbf{C}$ satisfying the following properties:
\begin{enumerate}
\item there exists $\Li{w}_{1}\in \RR_{+}$ such that for any $(x,y)\in E^{2}$, $u\in[0,\tstar{x}\wedge\tstar{y}]$, one has
$$\big| w(\ph{x}{u})-w(\ph{y}{u})\big| \leq \Li{w}_{1} |x-y|,$$
\item there exists $\Li{w}_{2}\in \RR_{+}$ such that for any $x\in E$, and $(t,s)\in [0,\tstar{x}]^{2}$, one has
$$\big| w(\ph{x}{t})-w(\ph{x}{s})\big| \leq \Li{w}_{2} |t-s|,$$
\item there exists $\Li{w}_{*}\in \RR_{+}$ such that for any $(x,y)\in E^{2}$, one has
$$\big| w(\ph{x}{\tstar{x}})-w(\ph{y}{\tstar{y}})\big| \leq \Li{w}_{*} |x-y|.$$
\end{enumerate}

In the sequel, for any function $w$ in $\mathbf{C}$, we denote by $C_{w}$ its bound:
\begin{equation*}
C_{w}=\sup_{x\in E}|w(x)|.
\end{equation*}
The following assumptions will be in force throughout.

\begin{assumption}
\label{H2}
The jump rate $\lambda$ is bounded and there exists $\Li{\lambda}_{1}\in \RR_{+}$ such that for any
$(x,y)\in E^{2}$, $u\in[0,\tstar{x}\wedge\tstar{y}[$,
$$\big| \lambda(\ph{x}{u})-\lambda(\ph{y}{u})\big| \leq \Li{\lambda}_{1} |x-y|.$$
\end{assumption}

\begin{assumption}
\label{H3}
The exit time $t^{*}$ is bounded and Lipschitz-continuous on $E$.
\end{assumption}

\begin{assumption}
\label{H5}
The Markov kernel $Q$ is Lipschitz in the following sense: there exists $\Li{Q}\in \RR_{+}$ such that for any function $w\in \mathbf{L}$ the following two conditions are satisfied:
\begin{enumerate}
\item for any $(x,y)\in E^{2}$, $u\in[0,\tstar{x}\wedge\tstar{y}]$, one has
$$\big| Qw(\ph{x}{u})-Qw(\ph{y}{u})\big| \leq \Li{Q} \Li{w}_{1} |x-y|,$$
\item for any $(x,y)\in E^{2}$, one has
$$\big| Qw(\ph{x}{\tstar{x}})-Qw(\ph{y}{\tstar{y}})\big| \leq \Li{Q} \Li{w}_{*} |x-y|.$$
\end{enumerate}
\end{assumption}

\section{Quantization}
\label{quanti}
The aim of this section is to describe the quantization procedure for a random variable and to recall some important properties that will be used in the sequel. There exists an extensive literature on quantization methods for random variables and processes. We do not pretend to present here an exhaustive panorama of these methods. However, the interested reader may for instance, consult the following works \cite{gray98, pages98, pages04} and references therein.
Consider $X$ an $\mathbb{R}^q$-valued random variable such that $\big\| X \big\|_p < \infty$ where $\big\| X \big\|_p$ denotes the $L_{p}$-nom of $X$: $\big\| X \big\|_p=\Big( \mathbb{E}[|X|^{p}]\Big)^{1/p}$. 

\bigskip

\noindent
Let $K$ be a fixed integer, the optimal $L_{p}$-quantization of the random variable $X$ consists in finding the best possible $L_{p}$-approximation of $X$ by a random vector $\widehat{X}$ taking at most $N$ values: $\widehat{X}\in \{x^{1},\ldots,x^{K}\}$.
This procedure consists in the following two steps:
\begin{enumerate}
\item Find a finite weighted grid $\Gamma\subset \mathbb{R}^q$ with $\Gamma= \{x^{1},\ldots,x^{K}\}$.
\item Set $\widehat{X}=\widehat{X}^{\Gamma}$ where $\widehat{X}^{\Gamma}=p_{\Gamma}(X)$ with $p_{\Gamma}$ denotes the closest neighbour projection on $\Gamma$.
\end{enumerate}

\noindent
The asymptotic properties of the $L_{p}$-quantization are given by the following result, see e.g. \cite{pages98}.
\begin{theorem}
If $\mathbb{E}[|X|^{p+\eta}]<+\infty$ for some $\eta>0$ then one has
\begin{eqnarray*}
\lim_{K\rightarrow \infty} K^{p/q} \min_{|\Gamma|\leq K} \| X-\widehat{X}^{\Gamma}\|^{p}_{p}& =& J_{p,q} 
\int |h|^{q/(q+p)}(u)du,
\end{eqnarray*}
where the law of $X$ is $P_{X}(du)=h(u) \lambda_{q}(du)+\nu$ with $\nu\perp\lambda_{d}$, $J_{p,d}$ a constant and $ \lambda_{q}$ the Lebesgue measure in $\mathbb{R}^{q}$.
\end{theorem}
Remark that $X$ needs to have finite moments up to the order $p+\eta$ to ensure the above convergence.
There exists a similar procedure for the optimal quantization of a Markov chain $\{X_{k}\}_{k\in \NN}$. There are two approaches to provide the quantized approximation of a Markov chain. The first one, based on the quantization at each time $k$ of the random variable $X_{k}$ is called the \textit{marginal quantization}. The second one that enhances the preservation of the Markov property is called \textit{Markovian quantization}. Remark that for the latter, the quantized Markov process is not homogeneous. These two methods are described in details in \cite[section 3]{pages04}. In this work, we used the marginal quantization approach for simplicity reasons.
\section{Impulse control problem}
\label{section impulse}
The formal probabilistic apparatus necessary to precisely define the impulse control problem is rather cumbersome, and will not be used in the sequel, therefore, for the sake of simplicity, we only present a rough description of the problem. 
The interested reader is referred to \cite{CD89} for a rigorous definition.

\bigskip

A strategy $\mathcal{S}=(\tau_n,R_n)_{n\geq 1}$ is a sequence of non-anticipative  intervention times $(\tau_n)_{n\geq 1}$ and non-anticipative $E$-valued random variables $(R_n)_{n\geq 1}$ on a measurable space $(\overline{\Omega},\overline{\mathcal{F}})$. 
Between the intervention times $\tau_i$ and $\tau_{i+1}$, the motion of the system is determined by the PDMP $\{X(t)\}$ starting from $R_i$.
If an intervention takes place at $x\in E$, then the set of admissible points where the decision-maker can send the system to is denoted by $\mathbb{U}\subset E$. We suppose that the control set $\mathbb{U}$ is finite and does not depend on $x$. The cardinal of the set $\mathbb{U}$ is denoted by $u$: 
\begin{equation*}
\mathbb{U} = \big\{ y^{i}: 1\leq i \leq u \big\}.
\end{equation*}
The strategy $\mathcal{S}$ induces a family of probability measures $P^{\mathcal{S}}_x$, $x\in E$, on $(\overline{\Omega},\overline{\mathcal{F}})$.
We define the class $\mathbb{S}$ of admissible strategies as the strategies $\mathcal{S}$ which satisfy $\tau_{\infty}=\infty$ $P^{\mathcal{S}}_x$-a.s. for all $x\in E$.

\bigskip

Associated to the strategy $\mathcal{S}$, we define the following discounted cost for a process starting at $x\in E$
\begin{equation*}
\mathcal{J}^{\mathcal{S}}(x)=E_x^{\mathcal{S}}\left[\int_0^{\infty}\e^{-\alpha s}f(Y_s)ds+\sum_{i=1}^{\infty}\e^{-\alpha \tau_i}c(Y_{\tau_i},Y_{\tau_i^+})\right],
\end{equation*}
where $E_x^{\mathcal{S}}$ is the expectation with respect to $P^{\mathcal{S}}_x$ and $\{Y_t\}$ is the process with interventions. The function $f$ then corresponds to the running cost and $c(x,y)$ corresponds to the intervention cost of moving the process from $x$ to $y$, $\alpha$ is a positive discount factor. We make the following assumption on the cost functions.

\begin{assumption}
\label{H6}
$f$ is a positive function in $\mathbf{L}$.
\end{assumption} 

\begin{assumption}
\label{H7}
The function $c$ is continuous on $\overline{E}\times\mathbb{U}$ and there exist $\Li{c}_{1}\in \RR_{+}$, $\Li{c}_{2}\in \RR_{+}$ and $\Li{c}_{*}\in \RR_{+}$ such that
\begin{enumerate}
\item for any $(x,y)\in E^{2}$, $u\in[0,\tstar{x}\wedge\tstar{y}]$, 
$$\max_{z\in \mathbb{U}} \big| c(\ph{x}{u},z)-c(\ph{y}{u},z)\big| \leq \Li{c}_{1} |x-y|,$$
\item for any $x\in E$, and $(t,s)\in [0,\tstar{x}]^{2}$, 
$$\max_{z\in \mathbb{U}} \big| c(\ph{x}{t},z)-c(\ph{x}{s},z)\big| \leq \Li{c}_{2} |t-s|,$$
\item for any $(x,y)\in E^{2}$, 
$$ \max_{z\in \mathbb{U}} \big| c(\ph{x}{\tstar{x}},z)-c(\ph{y}{\tstar{y}},z)\big| \leq \Li{c}_{*} |x-y|,$$
\item for any $(x,y)\in \overline{E}\times\mathbb{U}$, $0<c_0\leq c(x,y)\leq C_c$,
\item for any $(x,y,z)\in \overline{E}\times\mathbb{U}\times\mathbb{U}$, 
$$c(x,y)+c(y,z)\geq c(x,z).$$
\end{enumerate}
\end{assumption} 
The last assumption implies that the cost of taking two or more interventions instantaneously will not be lower than taking a single intervention. Finally, the value function for the discounted infinite horizon impulse control problem is defined for all $x$ in $E$ by
\begin{equation*}
\mathcal{V}(x)=\inf_{\mathcal{S}\in\mathbb{S}}\mathcal{J}^{\mathcal{S}}(x).
\end{equation*}
Associated to this impulse control problem, we define the following operators.
For $x\in E$, $t\geq 0$, $(v,w)\in \mathbf{C}^{2}$, set
\begin{eqnarray*}
F(x,t)&=&\int_{0}^{t\wedge \tstar{x}}\expla{x}{s}f\big(\ph{x}{s}\big)ds,\\
Hv(x,t)&=&\expla{x}{t\wedge \tstar{x}}v\big(\ph{x}{t\wedge \tstar{x}}\big)\\
& = & \mathbb{E}\Big[\e^{-\alpha (t\wedge \tstar{Z_{0}})} v\big(\ph{Z_{0}}{t\wedge \tstar{Z_{0}}}\big) 
\1_{\{S_{1}\geq t\wedge \tstar{Z_{0}}\}} \Big| Z_{0}=x \Big], \\
Iw(x,t)&=& \int_{0}^{t\wedge \tstar{x}}\expla{x}{s}\lambda Qw\big(\ph{x}{s}\big)ds, \\
& = & \mathbb{E}\Big[\e^{-\alpha S_{1}} w(Z_{1})  \1_{\{S_{1}< t\wedge \tstar{Z_{0}}\}} \Big| Z_{0}=x \Big].
\end{eqnarray*}
Finally for notational convenience, let us introduce for $(v,w) \in \mathbf{C}^{2}$, $x\in E$ and $t\geq0$.
\begin{eqnarray*}
J(v,w)(x,t)&=&F(x,t)+Hv(x,t)+Iw(x,t),\label{defJ}\\
Kw(x)&=&F(x,\tstar{x})+HQw(x,\tstar{x})+Iw(x,\tstar{x}),\label{defK}
\end{eqnarray*}
It is easy to show that for all $n\in \NN$
\begin{eqnarray}
Kv(x) & = & \mathbb{E} \Big[ F(Z_{n},\tstar{Z_{n}})+ \e^{-\alpha S_{n+1}} v(Z_{n+1})\big| Z_{n}= x\Big],
\label{DefKeq}\\
J(v,w)(x,t) & = & \mathbb{E} \Big[ F(Z_{n}, t)+ \e^{-\alpha S_{n+1}} w(Z_{n+1})\1_{\{S_{n+1}< t\wedge \tstar{Z_{n}}\}}\nonumber\\
&&+\e^{-\alpha t\wedge\tstar{Z_{n}}} v(\ph{Z_n}{t\wedge\tstar{Z_{n}}}\1_{\{S_{n+1}\geq t\wedge \tstar{Z_{n}}\}}
\big| Z_{n}= x\Big].
\label{DefJeq}
\end{eqnarray}
Note that these operators involve the original non controlled process $\{X(t)\}$ and only depend on the underlying Markov chain $(\Theta_n)=(Z_n,S_n)$. The equalities above are valid for all $n$ because $(\Theta_n)$ is an homogeneous Markov chain. Finally, for $(v,w) \in \mathbf{C}^{2}$, $\varphi$ defined on $\mathbb{U}$ and $x\in E$, set
\begin{eqnarray*}
M\varphi(x)& = & \inf_{y\in\mathbb{U}} \big\{ c(x,y)+\varphi(y) \big\}, \nonumber\\
L(v,w)(x) & = & \inf_{t\in\RR_{+}}J(v,w)(x,t) \wedge Kw(x), \label{DefL}\\
\mathcal{L}w(x)&=&L(Mw,w)(x). \nonumber
\end{eqnarray*}
As explained in \cite{CD89}, operator $\mathcal{L}$ applied to $w$ is the value function of the single-jump-or-intervention problem with cost function $w$ and the value function $\mathcal{V}$ can be computed by iterating  $\mathcal{L}$. 
More precisely, let $h$ be the cost associated to the no-impulse strategy: 
\begin{equation*}
h(x)=E_x\left[\int_0^{\infty}\e^{-\alpha s}f(X_s)ds\right],
\end{equation*}
for all $x\in E$.
Then we recall proposition~4 of \cite{CD89}.
\begin{proposition}\label{prop CD}
Assume that $g$ is in $\mathbf{L}$ and $g\geq h$. Define $\mathcal{V}_{0}^g=g$ and $\mathcal{V}_{n+1}^g=\mathcal{L}(\mathcal{V}_n^g)$, for all $n\geq 0$. Then for all $x\in E$
\begin{equation*}
\mathcal{V}(x)=\lim_{n\rightarrow\infty}\mathcal{V}_n^g(x).
\end{equation*}
\end{proposition}
As pointed out in \cite{CD89}, if one chooses exactly $g= h$, then $\mathcal{V}_n^h$ corresponds to the value function of the impulse problem where only $n$ jumps plus interventions are allowed, and after that, there are no further interventions.

\begin{remark}
Note that operator $L$ is quite similar to the operator used in optimal stopping, see e.g. \cite{costa88,AAP10}. However, the iteration procedure here does not rely on $L$ but on $\mathcal{L}$. The difference between operators $L$ and $\mathcal{L}$ comes from the operator $M$ that chooses optimally the next starting point. This is one of the main technical differences between approximating the value functions of an optimal stopping and impulse problems, and it makes the approximation scheme significantly more difficult, as explained in the next section.
\end{remark}

\section{Approximation of the value function}
\label{ApproxValue}
From now on, we assume that the distribution of $X(0)$ is given by $\delta_{x_0}$ for some fixed point $x_0$ in the state space $E$. We also choose a function $g$ in $\mathbf{L}$ satisfying $g\geq h$. Our approximation of the value function at $x_0$ is based on Proposition~\ref{prop CD}. Following the approach proposed by M. Davis and O. Costa in \cite{CD89}, we suppose now that we have selected a suitable index $N$ such that $\mathcal{V}(x_0)-\mathcal{V}_N^g(x_0)$ is small enough see the example in section \ref{section appli}. We turn to the approximation of $\mathcal{V}_N^g(x_0)$ which is the main object of this paper. In all generality, finding an index $N$ such that $\mathcal{V}(x_0)-\mathcal{V}_N^g(x_0)$ is below a prescribed level is a very difficult problem to solve. However, in particular cases one can hope to be able to evaluate the distance between $\mathcal{V}(x_0)$ and $\mathcal{V}_N^g(x_0)$.
As suggested by M. Davis and O. Costa in  \cite{CD89},
a value of $N$ can be chosen by calculating $\mathcal{V}_n^g(x_0)$ for different values of $n$ and stopping when the difference between two consecutive values is small enough. 
Our results of convergence are derived for a fixed but arbitrary $N$.


\bigskip

 Recall that if $\mathcal{V}_0= h$, then $\mathcal{V}_N^h$ corresponds to the value function of the impulse problem where only $N$ jumps plus interventions are allowed. This is an interesting problem to be solved in itself. For notational convenience, we will change our notation in the sequel and reverse the indices for the sequence $(\mathcal{V}_n^g)_{0\leq n\leq N}$. Set 
\begin{equation*}
\left\{\begin{array}{lcl}
v_N&=&g\ =\ \mathcal{V}_{0}^g, \\
v_n&=&\mathcal{L}v_{n+1} \ =\ \mathcal{V}_{N-n}^g,\quad\textrm{ for all }0\leq n<N.
\end{array}\right.
\end{equation*}

As explained in the introduction, the keystone of the approximation procedure for optimal stopping in \cite{AAP10} is that the analogue of Proposition~\ref{prop CD} yields a recursive construction of the random variables ${v}_n(Z_n)$.
Unfortunately, this key and important property does not hold anymore here. Indeed, one has:
\begin{eqnarray*}
{v}_{n}(Z_{n})&=&\mathcal{L}{v}_{n+1}(Z_{n})\\
&=&\bigg(\inf_{t\in\mathbb{R}_+}\mathbb{E} \Big[ F(Z_{n}, t)+ \e^{-\alpha S_{n+1}} v_{n+1}(Z_{n+1})\1_{\{S_{n+1}< t\wedge \tstar{Z_{n}}\}}\\
&&+\e^{-\alpha t\wedge\tstar{Z_{n}}} Mv_{n+1}(\ph{Z_n}{t\wedge\tstar{Z_{n}}}\1_{\{S_{n+1}\geq t\wedge \tstar{Z_{n}}\}}
\big| Z_{n}\Big]\bigg)\\
&&\wedge\mathbb{E} \Big[ F(Z_{n},\tstar{Z_{n}})+ \e^{-\alpha S_{n+1}} v_{n+1}(Z_{n+1})\big| Z_{n}\Big].
\end{eqnarray*}
And $Mv_{n+1}(\ph{Z_n}{t\wedge\tstar{Z_{n}}}$ cannot be written as a function of $v_{n+1}(Z_{n+1})$. Hence, we have no recursive construction of the random variables ${v}_n(Z_n)$ and we cannot apply the same procedure that we used for optimal stopping. Thus, we propose a new procedure to evaluate $Mv_{n+1}(\ph{Z_n}{t\wedge\tstar{Z_{n}}}$ separately from the main computation of the value function.

\bigskip

Note that for all $0\leq n<N$, to compute $Mv_{n+1}$ at any point, one actually only needs to evaluate the value functions $v_{n+1}$ at the points of the control grid $\mathbb{U}$. We propose again a recursive computation based on the Markov chain $(Z_n,S_n)$ but with a different starting point. Set $Z_0^y=y\in \mathbb{U}$ and $S_0^y=0$. We denote by $(Z_n^y,S_n^y)$ the Markov chain starting from this point $(y,0)$. One clearly knows $v_N=g$ on $\mathbb{U}$. Now suppose we have computed all the $v_n$ on $\mathbb{U}$ for $k+1\leq n\leq N$. Therefore, all functions $Mv_{n}$ are known everywhere. We can then propose the following recursive computation to evaluate $v_k$ at $y\in\mathbb{U}$:
\begin{equation}\label{rec Mv}
\left\{\begin{array}{l}
v_N(Z_{N-k}^y)= g(Z_{N-k}^y)\\
v_{k+n}(Z_n^y)=\mathcal{L}(Mv_{k+n+1},v_{k+n+1})(Z_n^y),\quad\textrm{ for all }0\leq n\leq N-k-1.   
\end{array}\right.
\end{equation}
This way, one obtains $v_k(Z_0^y)$ that exactly equals $v_k(y)$. Note that, since the functions $Mv_{k+n}$ are known, this provides a tractable recurrence relation on the random variables $v_{k+n}(Z_k^y)$.

\begin{remark}\label{remark g}
Note that this procedure requires the knowledge of function $g$ for all the random variables $(Z_n^y)_{n\leq N-1}$ defined for the different starting points $y\in\mathbb{U}$. This is why, in general, we are not able to use the no-impulse cost function $h$. Indeed, it is hard to compute this function, especially if we need to know it everywhere on the state space. The most practical solution is to take $g$ equal to a upper bound of $h$, and therefore constant.
\end{remark}

There is yet another new difficulty hidden in the recurrence relation (\ref{rec Mv}) above as regards its discretization. Indeed, to compute $v_n(y)$, one needs first to compute all the $v_{k+n}(Z_n^y)$ with $1\leq n\leq N-k$, and to compute $v_{k+1}(y)$ for instance, one has already computed all the $v_{k+n}(Z_{n-1}^y)$ for $2\leq n\leq N-k$. Unfortunately, one cannot re-use the values of $v_{k+n}(Z_{n-1}^y)$ to compute that of $v_{k+n}(Z_n^y)$, so the computation has to be started all over again each time, and one has to be very careful in the design of the approximation scheme. However, all these computations can be done with the same discretization grids for $(Z_n^y,S_n^y)$, so that our procedure is still reasonably fast, see section~\ref{section control grid} for details, and figure \ref{algo} for a graphical illustration of our procedure.

\begin{remark}
The recursive procedure (\ref{rec Mv}) is triangular in the sense that one needs to compute all the $v_{k+n}(Z_n^y)$ for $0\leq k\leq N$ and $0\leq n\leq N-k$.
\end{remark}

Our approximation procedure is in three steps, as explained in the following sections. The first step consists in replacing the continuous minimization in the definition of operator $L$ by a discrete-time minimization, on path adapted grids. The second step is specific to the impulse problem, and is due to the operator $M$ as explained in details above. The second step hence consists in carefully approximating the value functions $v_n$ on the control grid $\mathbb{U}$. The last step will then be similar to the approximation of the optimal stopping problem and will consist in approximating the value functions at the points of the quantization grids of the \emph{no impulse} process.

\subsection{Time discretization}
We define the path-adapted discretization grids as follows.

\begin{definition}
\label{grid}
For $z\in E$, set $\Delta(z)\in ]0,\tstar{z}[$. Define $ n(z)=\text{int}\big(\frac{\tstar{z}}{\Delta(z)}\big)-1$, where $\text{int}(x)$ denotes the greatest integer smaller than or equal to $x$.
The set of points $(t_{i})_{i\in\{0,\ldots,n(z)\}}$ with $\ds t_{i} = i\Delta(z)$ is denoted by $G(z)$.
This is the grid associated to the time interval $[0,\tstar{z}]$.
\end{definition}

\begin{remark}
\label{rq grille tstar}
It is important to note that, for all $z\in E$, not only one has $\tstar{z} \notin G(z)$, but also $\max G(z)=t_{n(z)}\leq \tstar{z}-\Delta(z)$. This property is crucial for the sequel.
\end{remark}
We propose the following approximation of operator $L$, where the continuous minimization is replaced by a discrete-time minimization on the path-adapted grids.

\begin{definition}
For $(v,w)\in \mathbf{L}^{2}$ and  $x\in E$, set
\begin{equation*}
L^{d}(v,w)(x)  = \min_{t \in G(x)} J(v,w)(x,t) \wedge Kw(x).
\end{equation*}
\end{definition}
Now we compute the error induced by the replacement of the continuous minimization by the discrete one.

\begin{lemma}\label{inf min J}
Let $(v,w)\in \mathbf{L}^{2}$. Then for all $x\in E$,
\begin{equation*}
\big|\inf_{t\leq t^*(x)}J(v,w)(x,t)-\min_{s\in G(x)}J(v,w)(x,s)\big| \leq\Big( C_f + C_wC_{\lambda} + \Li{v}_2+C_v(C_{\lambda}+\alpha)\Big)
\Delta(x).
\end{equation*}
\end{lemma}

\noindent{\textbf{Proof:}} 
We have
\begin{eqnarray*}
\big|\inf_{t\leq t^*(x)}J(v,w)(x,t)-\min_{s\in G(x)}J(v,w)(x,s)\big| & = &
\min_{s\in G(x)}J(v,w)(x,s)-\inf_{t\leq t^*(x)}J(v,w)(x,t).
\end{eqnarray*}
Clearly, there exists $\widebar{t} \in [0,\tstar{x}]$ such that $\ds \inf_{t\leq t^*(x)}J(v,w)(x,t)=J(v,w)(x,\widebar{t})$.
Moreover, there exists $0\leq i\leq n(x)$ such that $\widebar{t} \in [t_i,t_{i+1}]$ (with $t_{n(x)+1}=\tstar{x}$).
Consequently, Lemma~\ref{TechJt} yields
\begin{align*}
\big|\inf_{t\leq t^*(x)}J(v,w)(x,t)-\min_{s\in G(x)}J(v,w) & (x,s)\big| \leq J(v,w)(x,t_i)-J(v,w)(x,\widebar{t})\\
&\leq \Big( C_f + C_wC_{\lambda} + \Li{v}_2+C_v(C_{\lambda}+\alpha)\Big) |\widebar{t}-t_i|.
\end{align*}
implying the result.
\hfill $\Box$

\begin{lemma}
\label{lem erreur 2}
Let $(v,w)\in \mathbf{L}^{2}$ be nonnegative functions. Then for all $x\in E$,
\begin{equation*}
\big| L(v,w)(x)-L^{d}(v,w)(x) \big| \leq \Big( C_f + C_wC_{\lambda} + \Li{v}_2+C_v(C_{\lambda}+\alpha)\Big)
\Delta(x).
\end{equation*}
\end{lemma}

\noindent{\textbf{Proof:}} 
Since the functions $v$ and $w$ are nonnegative, it follows from the definition of $L$ and $L^{d}$ that
\begin{align*}
\big| L(v,w) & (x)-L^{d}(v,w)(x)\big| \leq \big| \inf_{t\leq t^*(x)}J(v,w)(x,t)-\min_{s\in G(x)}J(v,w)(x,s)\big|.
\end{align*}
Now in view of the previous lemma, one obtains the result.
\hfill $\Box$
\subsection{Approximation of the value functions on the control grid $\mathbb{U}$}
\label{section control grid}
We now need to introduce the quantized approximations of the underlying Markov chains $(\Theta_n^y)$. More precisely, we need several approximations at this stage, one for each starting point $y$ in the control set $\mathbb{U}$. 
Recall that $\mathbb{U}=\{y^i,\ 1\leq i\leq u\}$. For all $1\leq i\leq u$, let $({Z}_n^i,{S}_n^i)_{0\leq n\leq N-1}$ be the Markov chain $({Z}_n,{S}_n)_{0\leq n\leq N-1}$ with starting point $Z_0=y^i$, $S_0=0$, and let $(\widehat{Z}^i_n,\widehat{S}^i_n)_{0\leq n\leq N-1}$ be the quantized approximation of the sequence $({Z}_n^i,{S}_n^i)_{0\leq n\leq N-1}$, see Section~\ref{quanti}. The quantization algorithm provides us with a finite grid $\Gamma^{i,\Theta}_n\subset E\times \mathbb{R}_+$ at each time $0\leq n\leq N-1$ as well as weights for each point of the grid and transition probabilities from one grid to the next one, see e.g. \cite{bally03, pages98, pages04} for details. Set $p\geq 1$ such that $\Theta_n$ has finite moments at least up to the order $p+\epsilon$ for some positive $\epsilon$  and let $p^i_n$ be the closest-neighbour projection from $E\times \mathbb{R}_+$ onto $\Gamma^{i,\Theta}_n$ (for the distance of norm $p$; if there are several equally close neighbours, pick the one with the smallest index). Then the quantization of $\Theta_n^i$ conditionally to $Z_0=y^i$ is defined by
\begin{equation*}
\widehat{\Theta}^i_{n}=\big(\widehat{Z}^i_{n},\widehat{S}^i_{n}\big)=p^i_n\big({Z}_{n}^i,{S}_{n}^i\big).
\end{equation*}
We will also denote $\Gamma^{i,Z}_{n}$ the projection of $\Gamma^{i,\Theta}_n$ on $E$ and $\Gamma^{i,S}_{n}$ the projection of $\Gamma^{i,\Theta}_n$ on $\mathbb{R}_+$. 

\bigskip

Although $({Z}_n^i,{S}_n^i)$ is a Markov chain, its quantized approximation is usually not a Markov chain. It can be turned into a Markov chain by slightly changing the ponderations in the grids, see \cite{pages04b}, but this Markov chain will not be homogeneous in any case. Therefore, the following quantized approximations of operators $H$, $I$, $K$, $J$ and $L^d$ depend on both indices $n$ and $i$.
\begin{definition}
For $v\in \mathbf{L}^{2}$, $w$ defined on $\Gamma^{i,Z}_{n+1}$, $x\in E$, $0\leq n\leq N-1$, $1\leq i\leq u$ and $z\in \Gamma^{i,Z}_{n}$, consider 
\begin{align*}
\h{H}^i_{n+1}v(z,t)& =  \mathbb{E}\Big[\e^{-\alpha (t\wedge \tstar{\h{Z}^i_{n}})} v\big(\ph{\h{Z}^i_{n}}{t\wedge \tstar{\h{Z}^i_{n}}}\big) 
\1_{\{\h{S}^i_{n+1}\geq t\wedge \tstar{\h{Z}^i_{n}}\}} \Big| \h{Z}^i_{n}=z \Big], \\
\h{I}^i_{n+1}w(z,t)& = \mathbb{E}\Big[\e^{-\alpha \h{S}^i_{n+1}} w(\h{Z}^i_{n+1})  \1_{\{\h{S}^i_{n+1}< t\wedge \tstar{\h{Z}^i_{n}}\}} \Big| \h{Z}^i_{n}=z \Big], \\
\h{K}^i_{n+1}w(z) & =  \mathbb{E} \Big[ F(\h{Z}^i_{n},\tstar{\h{Z}^i_{n}})+ \e^{-\alpha \h{S}^i_{n+1}} w(\h{Z}^i_{n+1})\big| \h{Z}^i_{n}= z\Big], \\
\h{J}^i_{n+1}(v,w)(z,t) & = \mathbb{E} \Big[ F(\h{Z}^i_{n},t)+ \e^{-\alpha \h{S}^i_{n+1}} w(\h{Z}^i_{n+1})
\1_{\{\h{S}^i_{n+1}<t\wedge \tstar{\h{Z}^i_{n}}\}} \big| \h{Z}^i_{n}= z\Big] \nonumber \\
& \phantom{=} + \mathbb{E} \Big[ \e^{-\alpha (t\wedge \tstar{\h{Z}^i_{n}})} v( \ph{\h{Z}^i_{n}}{t\wedge \tstar{\h{Z}^i_{n}} } ) 
\1_{\{\h{S}^i_{n+1} \geq t\wedge \tstar{\h{Z}^i_{n}}\}} \big| \h{Z}^i_{n}= z\Big], \\
\h{L}^{i,d}_{n+1}(v,w)(z) & = \min_{t \in G(z)} \h{J}^i_{n+1}(v,w)(z,t) \wedge \h{K}^i_{n+1}w(z).
\end{align*}
\end{definition}
Our approximation scheme goes backwards in time, in as much as it is initialized with computing $v_N$ at the points of the last quantization grids $\Gamma^{i,Z}_{N}$, then $v_{N-1}$ is computed on $\Gamma^{i,Z}_{N-1}$ and so on. 
\begin{definition}
\label{vchapi}
Set $\widetilde{v}_{N}(y^i)=g(y^i)$ for $1\leq i \leq u$. Then, for $1\leq k\leq N-1$ and $1\leq i \leq u$, set $\widetilde{v}_{k}(y^i)=\widehat{v}^{i,k}_{k}(y^i)$, where
\begin{eqnarray*}
\widehat{v}^{i,k}_{N}(z) & = & g(z),\qquad z \in \Gamma^{i,Z}_{N-k},\\
\widehat{v}^{i,k}_{k+n-1}(z) & = & \widehat{L}_{n}^{i,d}(M\widetilde{v}_{k+n},\widehat{v}^{k}_{k+n})(z),\qquad z \in \Gamma^{i,Z}_{n-1},\qquad n \in\{1,\ldots,N-k\}.
\end{eqnarray*}
\end{definition}
See figure \ref{algo} for a graphical illustration of this numerical procedure.
\begin{remark}
Note the use of both $\widetilde{v}_{k+n}$ and $\widehat{v}^{k}_{k+n}$ in the scheme above. This is due to the fact that we have to reset all our calculations for each value function $\widetilde{v}_k$ and cannot use the calculations made for e.g. $\widetilde{v}_{k+1}$ because the value functions are evaluated at different points, and are approximated with different discrete operators. This is mostly because the quantized process $(\widehat{Z}^i_n,\widehat{S}^i_n)$ is not an homogeneous Markov chain.
\end{remark}
We can now state our first result on the convergence rate of this approximation.

\begin{theorem}\label{th erreur control}
For all $1\leq k\leq N-1$, $0\leq n \leq N-k-1$ and $1\leq i\leq d$, suppose that $\Delta(z)$ for $z\in \Gamma^{i,Z}_{n}$ is such that \begin{eqnarray*}
\sqrt{\frac{d^{4} \|Z_{n}^i-\h{Z}_{n}^i\|_{p} + d^{5} \big\| S_{n+1}^i- \h{S}_{n+1}^i \big\|_{p}}{d^{3}}}&<& \min_{z\in \Gamma^{i,Z}_{n}}\{\Delta(z)\}.
\end{eqnarray*}
Then we have
\begin{eqnarray*}
\lefteqn{\|v_{k+n}(Z_n^i)-\widehat{v}^{i,k}_{k+n}(\h{Z}_n^i)\|_p}\\
& \leq &
\big\|v_{k+n+1}(Z_{n+1}^i)-\h{v}_{k+n+1}^{i,k}(\h{Z}_{n+1}^i)\big\|_p +
\max_{y\in \mathbb{U}} \big| v_{k+n+1}(y)-\widetilde{v}_{k+n+1}(y)\big| \\
&&+ d^{1}_{k,n} \|Z_{n}^i-\h{Z}_{n}^i\|_{p} +2\Li{v_{k+n+1}} \big\|Z_{n+1}^i-\h{Z}_{n+1}^i\big\|_p
+C_f \big\| S_{n+1}^i- \h{S}_{n+1}^i \big\|_{p}  \\
&&+ d^{2}_{k,n} \big\| \Delta(\widehat{Z}_{n}^i) \big\|_p +
 +2 \sqrt{d^{3} \big( d^{4} \|Z_{n}^i-\h{Z}_{n}^i\|_{p} + d^{5} \big\| S_{n+1}^i- \h{S}_{n+1}^i \big\|_{p} \big)},
\end{eqnarray*}
with
\begin{eqnarray*}
d^{1}_{k,n} & = & \bigg\{ \Li{Q}  \Li{v_{k+n+1}}_*+2E_3\bigg\}\vee\bigg\{C_c(E_1+\alpha\Li{t^{*}})+ 2(\Li{c}_{1}+\Li{c}_{2}\Li{t^{*}})\bigg\} \\
&&+[v_{k+n}]+\Li{Q}  \Li{v_{k+n+1}}_1 \frac{C_{\lambda} }{\alpha}+\frac{C_{f}}{\alpha}(E_1+E_2) ,\\
d^{2}_{k,n} & = & C_f + C_{v_{k+n+1}}C_{\lambda} + \Li{c}_2+(C_c+C_{v_{k+n+1}})(C_{\lambda}+\alpha), \\
d^{3} & = & \big(\frac{2C_{f}}{\alpha}+C_{c} \big) C_{\lambda}, \\
d^{4} & = & \frac{C_f}{\alpha}(1+\Li{t^{*}})+C_c\Li{t^{*}}, \\
d^{5} & = & 2\Big(2\frac{C_f}{\alpha}+C_c\Big).
\end{eqnarray*}
\end{theorem}
\begin{remark}
\label{larochelle}
Recall that $v_N=\widehat{v}^{i,k}_{N}=\widetilde{v}_{N}=g$. Hence, one has
$$\| v_{N}(Z_{N-k}^i)-\widehat{v}^{i,k}_{N}(\h{Z}_{N-k}^i)\|_p\leq [g]\big\|\h{Z}^{i}_{N}-{Z}^{i}_{N}\big\|_p
\quad \text{ and } \quad \max_{y\in \mathbb{U}} \big| v_{N}(y)-\widetilde{v}_{N}(y)\big|=0.$$
In addition, the quantization error $\|\Theta^{i}_{n}-\widehat{\Theta}^{i}_{n}\|_p $ goes to zero as the number of points in the grids goes to infinity, see e.g. \cite{pages98}.
Therefore, according to Definition \ref{vchapi} and by using an induction procedure
$\max_{y\in \mathbb{U}} \big| v_{k}(y)-\widetilde{v}_{k}(y)\big|$
can be made arbitrarily small by an adequate choice of the discretization parameters.
From a theoretical point of view, the error can be calculated by iterating the result of Theorem \ref{th erreur control}.
However, this result is not presented here because it would lead to an intricate expression. From a numerical point of view, a computer  can easily estimate this error as shown in the example of section \ref{section appli}.
\end{remark}
The proof is going to be detailed in the following sections. We first split the error into four terms. For all $1\leq k\leq N-1$, $0\leq n \leq N-k-1$ and $1\leq i\leq d$, we have
\begin{eqnarray*}
\|v_{k+n}(Z_n^i)-\widehat{v}^{i,k}_{k+n}(\h{Z}_n^i)\|_p & \leq & \sum_{j=1}^{4} \Upsilon_{j}^i,
\end{eqnarray*}
where
\begin{eqnarray*}
\Upsilon_{1}^i & = & \|v_{k+n}(Z_{n}^i)-v_{k+n}(\widehat{Z}_{n}^i)\|_p ,\\
\Upsilon_{2}^i & = & \|L(Mv_{k+n+1},v_{k+n+1})(\widehat{Z}_{n}^i)-L^{d}(Mv_{k+n+1},v_{k+n+1})(\widehat{Z}_{n}^i)\|_p ,\\
\Upsilon_{3}^i & = & \|L^{d}(Mv_{k+n+1},v_{k+n+1})(\widehat{Z}_{n}^i)-\widehat{L}_{n+1}^{i,d}(Mv_{k+n+1},v_{k+n+1})(\widehat{Z}_{n}^i)\|_p ,\\
\Upsilon_{4}^i & = & \|\widehat{L}_{n+1}^{i,d}(Mv_{k+n+1},v_{k+n+1})(\widehat{Z}_{n}^i)-\widehat{L}_{n+1}^{i,d}(M\widetilde{v}_{k+n+1}, \widehat{v}_{k+n+1}^{i,k})(\widehat{Z}_{n}^i)\|_p.
\end{eqnarray*}
The first two terms are easy enough to handle thanks to Corollary~\ref{cor v lip} and lemma~\ref{lem erreur 2}.

\begin{lemma}\label{lem erreur 1 i}
A upper bound for $\Upsilon_{1}^i$ is
\begin{equation*}
\|v_{k+n}(Z_{n}^i)-v_{k+n}(\widehat{Z}_{n}^i)\|_p\leq [v_{k+n}] \|Z_{n}^i-\widehat{Z}_{n}^i\|_p.
\end{equation*}
\end{lemma}

\begin{lemma}\label{lem erreur 2 i}
A upper bound for $\Upsilon_{2}^i$ is
\begin{eqnarray*}
\lefteqn{\|L(Mv_{k+n+1},v_{k+n+1})(\widehat{Z}_{n}^i)-L^{d}(Mv_{k+n+1},v_{k+n+1})(\widehat{Z}_{n}^i)\|_p}\\
& \leq & \Big( C_f + C_{v_{k+n+1}}C_{\lambda} + \Li{c}_2+(C_c+C_{v_{k+n+1}})(C_{\lambda}+\alpha)\Big)
\big\| \Delta(\widehat{Z}_{n}^i)\big\|_p .
\end{eqnarray*}
\end{lemma}

The fourth term is also easy enough to deal with as it is a mere comparison of two finite weighted sums.
\begin{lemma}
A upper bound for $\Upsilon_{4}^i$ is
\begin{eqnarray*}
\lefteqn{ \|\widehat{L}_{n+1}^{i,d} (Mv_{k+n+1},v_{k+n+1})(\widehat{Z}_{n}^i)-\widehat{L}_{n+1}^{i,d}(M\widetilde{v}_{k+n+1},\widehat{v}_{k+n+1}^{i,k})(\widehat{Z}_{n}^i)\|_p}\\
 &\leq& \Li{v_{k+n+1}} \big\|Z_{n+1}^i-\h{Z}_{n+1}^i\big\|_p+\big\|v_{k+n+1}(Z_{n+1}^i)-\h{v}_{k+n+1}^{i,k}(\h{Z}_{n+1}^i)\big\|_p\\
&&+ \max_{y\in \mathbb{U}} \big| v_{k+n+1}(y)-\widetilde{v}_{k+n+1}(y)\big|.
\end{eqnarray*}
\end{lemma}

\noindent{\textbf{Proof:}} 
We clearly have
\begin{align*}
 \big\| \widehat{L}_{n+1}^{d} & (Mv_{k+n+1},v_{k+n+1})(\widehat{Z}_{n}^i)-\widehat{L}_{n+1}^{i,d}(M\widetilde{v}_{k+n+1},\widehat{v}_{k+n+1}^{i,k})(\widehat{Z}_{n}^i) \big\|_p \nonumber \\
& \leq 
\Big\| \max_{t \in G(\widehat{Z}_{n}^i)} \big| \h{J}_{n+1}^i(Mv_{k+n+1},v_{k+n+1})(\widehat{Z}_{n}^i,t)
-\h{J}_{n+1}^i(M\widetilde{v}_{k+n+1},\h{v}_{k+n+1}^{i,k})(\widehat{Z}_{n}^i,t) \big| \Big\|_p \nonumber \\
& \phantom{\leq} \vee \Big\| \h{K}_{n+1}^iv_{k+n+1}(\widehat{Z}_{n}^i) -  \h{K}_{n+1}^i\h{v}_{k+n+1}^{i,k}(\widehat{Z}_{n}^i)  \Big\|_p \nonumber \\
&\leq \big\|\mathbb{E}[v_{k+n+1}(\h{Z}_{n+1}^i)-\widehat{v}_{k+n+1}^{i,k}(\h{Z}_{n+1}^i)\big|\h{Z}_{n}^i]\big\|_p \nonumber \\
& \phantom{\leq}
+  \Big\| \mathbb{E} \Big[
Mv_{k+n+1}\big(\ph{\h{Z}_{n}^i}{t\wedge \tstar{\h{Z}_{n}^i}}\big)- M\widetilde{v}_{k+n+1} \big(\ph{\h{Z}_{n}^i}{t\wedge \tstar{\h{Z}_{n}^i}}\big)
\Big| \h{Z}_{n}^i \Big] \Big\|_p
 \nonumber \\
&\leq \big\|v_{k+n+1}(\h{Z}_{n+1}^i)-v_{k+n+1}({Z}_{n+1}^i)\big\|_p+\big\|v_{k+n+1}({Z}_{n+1}^i)-\widehat{v}_{k+n+1}^{i,k}(\h{Z}_{n+1}^i)\big\|_p
\nonumber \\
& \phantom{\leq} + \max_{y\in \mathbb{U}} \big| v_{k+n+1}(y)-\widetilde{v}_{k+n+1}(y)\big|,
\end{align*}
showing the result.
\hfill $\Box$

\bigskip

We now turn to the third term. This is the key step of the error evaluation, because on the one hand, this is where we deal with the indicator functions. The main idea is that although they are not continuous, we prove in the following two lemmas that the set where the discontinuity actually occurs is of small enough probability. This is also where our special choice of time discretization grids is crucial. On the other hand, we use here the specific properties of quantization.
\begin{lemma}
\label{lem indic1i}
For all $1\leq i\leq d$, $n\in \{0,\ldots,N-1\}$ and $\ds 0<\eta<\min_{z\in \Gamma^{i,Z}_{n}}\{\Delta(z)\}$,
\begin{equation*}
\big\|\1_{\tstar{Z_{n}^i}<\tstar{\h{Z}_{n}^i}-\eta}\big\|_p \leq \frac{[t^*] \|Z_{n}^i-\h{Z}_{n}^i\|_p}{\eta}.
\end{equation*}
\end{lemma}

\noindent{\textbf{Proof:}} 
By using the Chebychev's inequality, one clearly has
\begin{eqnarray*}
\mathbb{E}\Big[|\1_{\tstar{Z_{n}^i}<\tstar{\h{Z}_{n}^i}-\eta}|^p\Big] &=&
\mathbf{P} \Big( \tstar{Z_{n}^i}<\tstar{\h{Z}_{n}^i}-\eta \Big) \\
&\leq&\mathbf{P} \Big( \big| \tstar{Z_n^i}-\tstar{\h{Z}_n^i} \big|>\eta\ \Big) \leq\  \frac{[t^*]^p\|Z_{n}^i-\h{Z}_{n}^i\|_p^p}{\eta^p},
\end{eqnarray*}
showing the result.\hfill $\Box$

\begin{lemma}
\label{lem indic2i}
For all $1\leq i\leq d$, $n\in \{0,\ldots,N-1\}$ and $\ds 0<\eta<\min_{z\in \Gamma^{i,Z}_{n}}\{\Delta(z)\}$,
\begin{align*}
\Big\| \max_{s\in G(\widehat{Z}_{n}^i)} \mathbb{E}\big[
| & \1_{\{{S}_{n+1}^i<s\wedge \tstar{Z_{n}^i}\}}-\1_{\{\h{S}_{n+1}^i<s\wedge \tstar{\h{Z}_{n}^i}\}}| \big| \h{Z}_{n}^i\big]
\Big\|_p \nonumber \\
& \leq  \frac{2}{\eta}\|{S}_{n+1}^i-\h{S}_{n+1}^i\|_p + C_{\lambda}\eta + \frac{2[t^*]}{\eta} \|Z_{n}^i-\h{Z}_{n}^i\|_p.
\end{align*}
\end{lemma}

\noindent{\textbf{Proof:}} 
Set $\ds 0<\eta<\min_{z\in \Gamma^{i,Z}_{n}}\{\Delta(z)\}$ and $s\in G(\widehat{Z}_{n}^i)$.
By definition of the grid $G(\widehat{Z}_{n}^i)$ and $\eta$, one has $s+\eta<\tstar{\h{Z}_{n}^i}$, see Remark~\ref{rq grille tstar}.
Thus, the difference of indicator functions can be written as
\begin{align*}
\big| & \1_{\{{S}_{n+1}^i<s\wedge \tstar{Z_{n}^i}\}}-\1_{\{\h{S}_{n+1}^i<s\wedge \tstar{\h{Z}_{n}^i}\}}\big| \nonumber \\
& \leq \big|\1_{\{{S}_{n+1}^i<s\wedge \tstar{Z_{n}^i}\}}-\1_{\{\h{S}_{n+1}^i<s\wedge \tstar{\h{Z}_{n}^i}\}}\big| 
\Big[ \1_{\{ \tstar{Z_{n}^i}\leq  \tstar{\h{Z}_{n}^i}-\frac{\eta}{2} \}}
+\1_{\{ \tstar{Z_{n}^i} >  \tstar{\h{Z}_{n}^i}-\frac{\eta}{2} \}} \Big] \nonumber \\
& \leq  \1_{\{ \tstar{Z_{n}^i}\leq  \tstar{\h{Z}_{n}^i}-\frac{\eta}{2} \}}
+ \1_{\{ \tstar{Z_{n}^i} >  s+\frac{\eta}{2} \}} \big|\1_{\{{S}_{n+1}^i<s\}}-\1_{\{\h{S}_{n+1}^i<s\}}\big| \nonumber \\
& \leq  \1_{\{ \tstar{Z_{n}^i}\leq  \tstar{\h{Z}_{n}^i}-\frac{\eta}{2} \}}
+\1_{\{|{S}_{n+1}^i-\h{S}_{n+1}^i|>\frac{\eta}{2}\}} + \1_{\{ \tstar{Z_{n}^i} >  s+\frac{\eta}{2} \}} \1_{\{|{S}_{n+1}^i-s|\leq\frac{\eta}{2}\}}.
\end{align*}
This yields
\begin{align}
\Big\| & \max_{s\in G(\widehat{Z}_{n}^i)} \mathbb{E}\big[
| \1_{\{{S}_{n+1}^i<s\wedge \tstar{Z_{n}^i}\}}-\1_{\{\h{S}_{n+1}^i<s\wedge \tstar{\h{Z}_{n}^i}\}}| \big| \h{Z}_{n}^i\big] \Big\|_p
\leq \big\| \1_{\{ \tstar{Z_{n}^i}\leq  \tstar{\h{Z}_{n}^i}-\frac{\eta}{2} \}} \big\|_{p}
\nonumber \\ 
& + \big\|\1_{\{|{S}_{n+1}^i-\h{S}_{n+1}^i|>\frac{\eta}{2}\}}\big\|_p+\big\|\max_{s\in G(\widehat{Z}_{n}^i)}
\mathbb{E}\big[  \1_{\{ \tstar{Z_{n}^i} >  s+\frac{\eta}{2} \}} \1_{\{|{S}_{n+1}^i-s|\leq\frac{\eta}{2}\}} \big| \h{Z}_{n}^i\big]\big\|_p.
\label{indic1EQ1i}
\end{align}
On the one hand, Chebychev's inequality gives
\begin{equation}
\big\| \1_{\{|{S}_{n+1}^i-\h{S}_{n+1}^i|>\frac{\eta}{2}\}} \big\|_{p}^{p} =
\mathbf{P}( |{S}_{n+1}^i-\h{S}_{n+1}^i|>\frac{\eta}{2} ) \leq \frac{2^{p} \big\| {S}_{n+1}^i-\h{S}_{n+1}^i \big\|_{p}^{p}}{\eta^p}.
\label{indic1EQ2i}
\end{equation}
On the other hand, one has
\begin{align}
\mathbb{E}\big[\1_{\{ \tstar{Z_{n}^i} >  s+\frac{\eta}{2} \}} \1_{\{|{S}_{n+1}^i-s|\leq\frac{\eta}{2}\}} \big| \h{Z}_{n}^i\big]
&= \mathbb{E}\Big[\mathbb{E}\big[ \1_{\{ \tstar{Z_{n}^i} >  s+\frac{\eta}{2} \}} \1_{\{s-\frac{\eta}{2}\leq {S}_{n+1}^i\leq s+\frac{\eta}{2}\}} | Z_{n}^i\big] \Big| \h{Z}_{n}^i\Big]\nonumber\\
&= \mathbb{E}\Big[\1_{\{ \tstar{Z_{n}^i} >  s+\frac{\eta}{2} \}} \int_{s-\frac{\eta}{2}}^{s+\frac{\eta}{2}}\lambda(\phi(Z_{n}^i,u))du\Big| \h{Z}_{n}^i\Big]
\nonumber\\
&\leq \eta C_{\lambda}.
\label{indic1EQ3i}
\end{align}
Combining Lemma \ref{lem indic1i} and equations (\ref{indic1EQ1i})-(\ref{indic1EQ3i}), the result follows.
\hfill $\Box$

\bigskip

We now look up the error made in replacing $K$ by $\h{K}_{n+1}^i$. This is where we use the specific properties of quantization.
\begin{lemma}
\label{Er3ai}
For all $1\leq i \leq d$, $k\in\{1,\ldots,N-1\}$ and $n \in\{1,\ldots,N-k\}$, one has
\begin{align*}
\Big\| K & v_{k+n+1}(\h{Z}_{n}^i) - \h{K}_{n+1}^iv_{k+n+1}(\h{Z}_{n}^i)\Big\|_{p}  \leq 
C_f \|S_{n+1}^i-\h{S}_{n+1}^i\|_{p} + \Li{v_{k+n+1}} \|Z_{n+1}^i-\h{Z}_{n+1}^i\|_{p} \nonumber \\
& +\bigg\{\Li{Q} \Li{v_{k+n+1}}_1\frac{C_{\lambda}}{\alpha}+\Li{Q} \Li{v_{k+n+1}}_{*} +\frac{C_{f}}{\alpha}\Big(E_{1}+E_{2}\Big)
+ 2 E_{3} \bigg\} \|Z_{n}^i-\h{Z}_{n}^i\|_{p}.
\end{align*}
\end{lemma}

\noindent{\textbf{Proof:}} 
We have
\begin{eqnarray}
\lefteqn{\Big| Kv_{k+n+1}(\h{Z}_{n}^i) - \h{K}_{n+1}^iv_{k+n+1}(\h{Z}_{n}^i)\Big|}\label{EQ1Er3ai}\\
 & \leq &
\Big| Kv_{k+n+1}(\h{Z}_{n}^i) - \mathbb{E} \big[ Kv_{k+n+1}(Z_{n}^i) \big| \h{Z}_{n}^i\big] \Big| +\Big| \mathbb{E} \big[ Kv_{k+n+1}(Z_{n}^i) \big| \h{Z}_{n}^i\big]- \h{K}_{n+1}^iv_{k+n+1}(\h{Z}_{n}^i)\Big| \nonumber \\
& \leq & \mathbb{E} \Big[ \big| Kv_{k+n+1}(\h{Z}_{n}^i) - Kv_{k+n+1}(Z_{n}^i)  \big| \Big| \h{Z}_{n}^i\Big] +\Big| \mathbb{E} \big[ Kv_{k+n+1}(Z_{n}^i) \big| \h{Z}_{n}^i\big]- \h{K}_{n+1}^iv_{k+n+1}(\h{Z}_{n}^i)\Big|.\nonumber
\end{eqnarray}
By using the Lipschitz property of $K$ stated in Lemma \ref{TechK}, we obtain
\begin{align}
\Big\|  \mathbb{E} \Big[ \big| K & v_{k+n+1}(\h{Z}_{n}^i) - Kv_{k+n+1}(Z_{n}^i)  \big| \Big| \h{Z}_{n}^i\Big] \Big\|_{p} \label{EQ3Er3ai} \\
& \leq
\bigg\{\Li{Q} \Li{v_{k+n+1}}_1\frac{C_{\lambda}}{\alpha}+\Li{Q} \Li{v_{k+n+1}}_{*} +C_{v_{k+n+1}}\Big(E_{1}+E_{2}\Big)
+E_{3} \bigg\} \|Z_{n}^i-\h{Z}_{n}^i\|_{p}.\nonumber
\end{align}
Then, recall that by construction of the quantized process, one has $\big(\widehat{Z}_{n}^i,\widehat{S}_{n}^i\big)=p_n^i\big({Z}_{n}^i,{S}_{n}^i\big)$.
Hence we have the following crucial property: $\sigma\{\h{Z}_{n}^i\} \subset \sigma\{{Z}_{n}^i, S_n^i\}$.
By using the special structure of the PDMP $\{X(t)\}$, we also have 
$\sigma\{{Z}_{n}^i, S_n^i\} \subset \mathcal{F}_{T_{n}}$, so that one has $\sigma\{\h{Z}_{n}^i\} \subset \sigma\{Z_{n}^i\}$.
It now follows from the definition of $K$ given in equation (\ref{DefKeq}) that
\begin{align}
\Big| \mathbb{E} \big[ Kv_{k+n+1}(Z_{n}^i) \big| \h{Z}_{n}^i\big]- &\h{K}_{n+1}^i  v_{k+n+1}(\h{Z}_{n}^i)\Big| \leq 
\mathbb{E} \Big[ \big| F(Z_{n}^i,\tstar{Z_{n}^i}) -  F(\h{Z}_{n}^i,\tstar{\h{Z}_{n}^i}) \big| \Big| \h{Z}_{n}^i \Big]\nonumber  \\
& + \mathbb{E} \Big[ \big| \e^{-\alpha S_{n+1}^i} v_{k+n+1}(Z_{n+1}^i)- \e^{-\alpha \h{S}_{n+1}^i} v_{k+n+1}(\h{Z}_{n+1}^i) \big| \Big| \h{Z}_{n}^i \Big].\label{EQ2Er3ai}
\end{align}
From Lemma \ref{TechFI}, we readily obtain
\begin{align}
\Big\| \mathbb{E}  \Big[ \big| F(Z_{n}^i,\tstar{Z_{n}^i}) -  F(\h{Z}_{n}^i,\tstar{\h{Z}_{n}^i}) \big| \Big| \h{Z}_{n}^i \Big] \Big\|_{p} 
& \leq E_{3} \|Z_{n}^i-\h{Z}_{n}^i\|_{p},
\label{EQ4Er3ai}
\end{align}
and it is easy to show that
\begin{align}
\Big\| \mathbb{E} \Big[ \big| \e^{-\alpha S_{n+1}^i} & v_{k+n+1}(Z_{n+1}^i)- \e^{-\alpha \h{S}_{n+1}^i}  v_{k+n+1}(\h{Z}_{n+1}^i) \big| \Big| \h{Z}_{n}^i \Big] \Big\|_{p}
\nonumber \\
& \leq \Li{v_{k+n+1}}  \|Z_{n+1}^i-\h{Z}_{n+1}^i\|_{p} + \alpha C_{v_{k+n+1}} \|S_{n+1}^i-\h{S}_{n+1}^i\|_{p}.
\label{EQ5Er3ai}
\end{align}
Finally, recalling that $C_{v_{k+n}}\leq \frac{C_{f}}{\alpha}$ and combining equations
(\ref{EQ1Er3ai})-(\ref{EQ5Er3ai}) we obtain the expected result. \hfill $\Box$

\bigskip

We turn to the error made in replacing $J$ by $\h{J}_{n+1}^i$. Here we use the specific properties of quantization again, and the lemmas on indicator functions.
\begin{lemma}
\label{Er3bi}
For all $1\leq i \leq d$, $k\in\{1,\ldots,N-1\}$, $n \in\{1,\ldots,N-k\}$, and $\ds 0<\eta<\min_{z\in \Gamma^{i,Z}_{n}}\{\Delta(z)\}$, one has
\begin{align*}
\Big\| \max_{t\in G(\h{Z}_{n}^i)} & \big| J(Mv_{k+n+1},v_{k+n+1})(\h{Z}_{n}^i,t)  - \h{J}_{n+1}^i(Mv_{k+n+1},v_{k+n+1})(\h{Z}_{n}^i,t) \big| \Big\|_{p}\\
& \leq  \bigg\{\Li{Q}  \Li{v_{k+n+1}}_1 \frac{C_{\lambda} }{\alpha} + \frac{C_{f}}{\alpha}(E_1+E_2+\alpha\Li{t^{*}})
+C_c(E_1+\alpha\Li{t^{*}})\\
& \phantom{\leq  \bigg\{}
+ 2(\Li{c}_{1}+\Li{c}_{2}\Li{t^{*}})+ \frac{1}{\eta}\Big(\frac{C_f}{\alpha}(1+\Li{t^{*}})+C_c\Li{t^{*}}\Big)\bigg\}\|Z_{n}^i-\h{Z}_{n}^i\|_{p} \\
& \phantom{\leq} +\Big\{ \frac{2}{\eta}\Big(2\frac{C_f}{\alpha}+C_c\Big)+C_f \Big\}
\big\| S_{n+1}^i- \h{S}_{n+1}^i \big\|_{p}
\nonumber \\
&\phantom{\leq} + \Li{v_{k+n+1}} \big\| Z_{n+1}^i- \h{Z}_{n+1}^i \big\|_{p}
+\big(\frac{2C_{f}}{\alpha}+C_{c} \big) C_{\lambda}\eta.
\end{align*}
\end{lemma}

\noindent{\textbf{Proof:}} 
By definition of $J$, we have
\begin{eqnarray}
\lefteqn{\Big| J(Mv_{k+n+1},v_{k+n+1})(\h{Z}_{n}^i,t)  - \h{J}_{n+1}^i(Mv_{k+n+1},  v_{k+n+1})(\h{Z}_{n}^i,t) \Big|}\label{EQ1Er3bi}\\
&  \leq &
\Big| Iv_{k+n+1}(\h{Z}_{n}^i,t)  - \h{I}_{n+1}^iv_{k+n+1}(\h{Z}_{n}^i,t) \Big| 
+\Big| HMv_{k+n+1}(\h{Z}_{n}^i,t)  - \h{H}_{n+1}^iMv_{k+n+1}(\h{Z}_{n}^i,t) \Big|.\nonumber
\end{eqnarray}
For the first term on the right hand side of equation (\ref{EQ1Er3bi}), we proceed as  for $K$ in the preceding lemma
\begin{eqnarray*}
\Big| Iv_{k+n+1}(\h{Z}_{n}^i,t)  - \h{I}_{n+1}^iv_{k+n+1}(\h{Z}_{n}^i,t) \Big|
& \leq & \mathbb{E} \Big[ \big| Iv_{k+n+1}(\h{Z}_{n}^i,t) - Iv_{k+n+1}(Z_{n}^i,t)  \big| \Big| \h{Z}_{n}^i\Big]\\
&&+\Big| \mathbb{E} \big[ Iv_{k+n+1}(Z_{n}^i,t) \big| \h{Z}_{n}^i\big]- \h{I}_{n+1}^iv_{k+n+1}(\h{Z}_{n}^i,t)\Big|.
\end{eqnarray*}
On the one hand, it follows from Lemma \ref{TechFI} that
\begin{eqnarray*}
\lefteqn{\Big\| \max_{t\in G(\h{Z}_{n}^i)}  \mathbb{E} \Big[ \big| Iv_{k+n+1}(\h{Z}_{n}^i,t) - Iv_{k+n+1}(Z_{n}^i,t)  \big| \Big| \h{Z}_{n}^i\Big] \Big\|_{p}}\\
& \leq&  \bigg\{ \frac{1}{\alpha} \Big( \Li{Q}  \Li{v_{k+n+1}}_1 C_{\lambda} +C_{v_{k+n+1}}\Li{\lambda}_1 \big(1+C_{\lambda}C_{t^{*}}\big)  \Big) +C_{v_{k+n+1}}C_{\lambda} \Li{t^{*}} \bigg\} \|Z_{n}^i-\h{Z}_{n}^i\|_{p}.
\end{eqnarray*}
On the other hand, we use again the fact that $\sigma\{\h{Z}_{n}^i\} \subset \sigma\{Z_{n}^i\}$ to obtain
\begin{align*}
\Big| \mathbb{E} \big[ & Iv_{k+n+1}(Z_{n}^i,t) \big| \h{Z}_{n}^i\big]- \h{I}_{n+1}^i v_{k+n+1}(\h{Z}_{n}^i,t)\Big| \nonumber \\
& \leq \mathbb{E} \Big[ \1_{\{S_{n+1}^i< t\wedge \tstar{Z_{n}^i}\}}
\big| \e^{-\alpha S_{n+1}^i} v_{k+n+1}(Z_{n+1}^i)- \e^{-\alpha \h{S}_{n+1}^i} v_{k+n+1}(\h{Z}_{n+1}^i) \big| \Big| \h{Z}_{n}^i \Big]\nonumber \\
& \phantom{\leq} + \mathbb{E} \Big[ \e^{-\alpha \h{S}_{n+1}^i} v_{k+n+1}(\h{Z}_{n+1}^i)
\big|  \1_{\{S_{n+1}^i< t\wedge \tstar{Z_{n}^i}\}}-\1_{\{\h{S}_{n+1}^i< t\wedge \tstar{\h{Z}_{n}^i}\}}  \big| \Big| \h{Z}_{n}^i \Big].
\end{align*}
It remains to deal with the indicator function. Lemma \ref{lem indic2i} yields
\begin{eqnarray*}
\lefteqn{\Big\| \max_{t\in G(\h{Z}_{n}^i)}
\big| \mathbb{E} \big[ Iv_{k+n+1}(Z_{n}^i,t) \big| \h{Z}_{n}^i\big]-  \h{I}_{n+1}^i v_{k+n+1}(\h{Z}_{n}^i,t)\big| \Big\|_{p}}\\
& \leq & C_{v_{k+n+1}}C_{\lambda}\eta+(\alpha+\frac{2}{\eta}) C_{v_{k+n+1}} \big\| S_{n+1}^i- \h{S}_{n+1}^i \big\|_{p}
\\
&& + \Li{v_{k+n+1}} \big\| Z_{n+1}^i- \h{Z}_{n+1}^i\big\|_{p}
+\frac{2C_{v_{k+n+1}}\Li{t^{*}}}{\eta} \big\| Z_{n}^i- \h{Z}_{n}^i \big\|_{p}.
\end{eqnarray*}
By using the same arguments and lemmas~\ref{lem M} and \ref{TechFI}, we obtain similar results for the second term on the right hand side of equation (\ref{EQ1Er3bi}), namely
\begin{align*}
\Big\| \max_{t\in G(\h{Z}_{n}^i)} & \mathbb{E} \Big[ \big| HMv_{k+n+1}(\h{Z}_{n}^i,t) - HMv_{k+n+1}(Z_{n}^i,t)  \big| \Big| \h{Z}_{n}^i\Big] \Big\|_{p} \nonumber \\
& \leq  \bigg\{ \Li{c}_{1}+\Li{c}_{2}\Li{t^{*}}+(C_{v_{k+n+1}}+C_{c})\big(C_{t^{*}}\Li{\lambda}_1+(C_{\lambda}+\alpha) \Li{t^{*}}\big) \bigg\} \|Z_{n}^i-\h{Z}_{n}^i\|_{p},
\end{align*}
and
\begin{align*}
\Big\| \max_{t\in G(\h{Z}_{n}^i)}
\big| \mathbb{E} & \big[ HMv_{k+n+1}(Z_{n}^i,t) \big| \h{Z}_{n}^i\big]-  \h{H}_{n+1}^i Mv_{k+n+1}(\h{Z}_{n}^i,t)\big| \Big\|_{p} \nonumber \\
& \leq \bigg\{ \alpha \Li{t^{*}} (C_{v_{k+n+1}}+C_{c}) + \Li{c}_{1}+\Li{c}_{2}\Li{t^{*}}
+\frac{2[t^*]}{\eta} (C_{v_{k+n+1}}+C_{c}) \bigg\}
\|Z_{n}^i-\h{Z}_{n}^i\|_p \nonumber \\
& \phantom{=} +\frac{2}{\eta}(C_{v_{k+n+1}}+C_{c}) \|{S}_{n+1}^i-\h{S}_{n+1}^i\|_p
+ C_{\lambda}(C_{v_{k+n+1}}+C_{c}) \eta,
\end{align*}
showing the result.
\hfill $\Box$

\bigskip

We now add up the preceding results to one obtains the following upper bound for $\Upsilon_{3}^i$.

\begin{lemma}\label{lem erreur 3 i}
A upper bound for $\Upsilon_{3}^i$ is
\begin{eqnarray*}
\lefteqn{\|L^{d}(Mv_{k+n+1},v_{k+n+1})(\widehat{Z}_{n}^i)-\widehat{L}_{n+1}^{i,d}(Mv_{k+n+1},v_{k+n+1})(\widehat{Z}_{n}^i)\|_p}\\
& \leq &\|Z_{n}^i-\h{Z}_{n}^i\|_{p}\Bigg\{ \Li{Q}  \Li{v_{k+n+1}}_1 \frac{C_{\lambda} }{\alpha}+\frac{C_{f}}{\alpha}(E_1+E_2)+\bigg\{ \Li{Q}  \Li{v_{k+n+1}}_*\\
&&+2E_3\bigg\}\vee\bigg\{C_c(E_1+\alpha\Li{t^{*}})+ 2(\Li{c}_{1}+\Li{c}_{2}\Li{t^{*}})+\frac{1}{\eta}\Big(\frac{C_f}{\alpha}(1+\Li{t^{*}})+C_c\Li{t^{*}}\Big)\bigg\}\Bigg\}\\
&&+\big\| S_{n+1}^i- \h{S}_{n+1}^i \big\|_{p}\Big\{C_f+ \frac{2}{\eta}\Big(2\frac{C_f}{\alpha}+C_c\Big)\Big\}\\
&&+\Li{v_{k+n+1}} \big\| Z_{n+1}^i- \h{Z}_{n+1}^i \big\|_{p}+\big(\frac{2C_{f}}{\alpha}+C_{c} \big) C_{\lambda}\eta.
\end{eqnarray*}
\end{lemma}

\subsection{Approximation of the value function}
Now we have computed the value functions on the control grid, we turn to the actual approximation of $v_0$. As in the preceding section, we define the quantized approximation of the underlying Markov chain $(\Theta_n)$ starting from $(x_0,0)$, the actual starting point of the PDMP. Let $(\widehat{Z}_n,\widehat{S}_n)_{0\leq n\leq N-1}$ be the quantized approximation of the sequence $({Z}_n,{S}_n)_{0\leq n\leq N-1}$. The quantization algorithm provides us with another series of finite grids $\Gamma^{\Theta}_n\subset E\times \mathbb{R}_+$ for all $0\leq n\leq N-1$ as well as weights for each point of the grids and transition probabilities from one grid to the next one. Let $p_n$ be the closest-neighbor projection from $E\times \mathbb{R}_+$ onto $\Gamma^{\Theta}_n$. Then the quantization of $\Theta_n$ conditionally to $Z_0=x_0$ is defined by
\begin{equation*}
\widehat{\Theta}_{n}=\big(\widehat{Z}_{n},\widehat{S}_{n}\big)=p_n\big({Z}_{n},{S}_{n}\big).
\end{equation*}
We will also denote $\Gamma^{Z}_{n}$ the projection of $\Gamma^{\Theta}_n$ on $E$ and $\Gamma^{S}_{n}$ the projection of $\Gamma^{\Theta}_n$ on $\mathbb{R}_+$. 
We use yet again new quantized approximations of operators $H$, $I$, $K$, $J$ and $L^d$.
\begin{definition}
For $v\in \mathbf{L}^{2}$, $w$ defined on $\Gamma^{Z}_{n+1}$, $x\in E$, $n\in\{0,\ldots,N-1\}$ and $z\in \Gamma^{z}_{n}$, consider 
\begin{align*}
\h{H}_{n+1}v(z,t)& =  \mathbb{E}\Big[\e^{-\alpha (t\wedge \tstar{\h{Z}_{n}})} v\big(\ph{\h{Z}_{n}}{t\wedge \tstar{\h{Z}_{n}}}\big) 
\1_{\{\h{S}_{n+1}\geq t\wedge \tstar{\h{Z}_{n}}\}} \Big| \h{Z}_{n}=z \Big], \\
\h{I}_{n+1}w(z,t)& = \mathbb{E}\Big[\e^{-\alpha \h{S}_{n+1}} w(\h{Z}_{n+1})  \1_{\{\h{S}_{n+1}< t\wedge \tstar{\h{Z}_{n}}\}} \Big| \h{Z}_{n}=z \Big] ,\\
\h{K}_{n+1}v(z) & =  \mathbb{E} \Big[ F(\h{Z}_{n},\tstar{\h{Z}_{n}})+ \e^{-\alpha \h{S}_{n+1}} v(\h{Z}_{n+1})\big| \h{Z}_{n}= z\Big] ,\\
\h{J}_{n+1}(v,w)(z,t) & = \mathbb{E} \Big[ F(\h{Z}_{n},t)+ \e^{-\alpha \h{S}_{n+1}} w(\h{Z}_{n+1})
\1_{\{\h{S}_{n+1}<t\wedge \tstar{\h{Z}_{n}}\}} \big| \h{Z}_{n}= z\Big] \nonumber \\
& \phantom{=} + \mathbb{E} \Big[ \e^{-\alpha (t\wedge \tstar{\h{Z}_{n}})} v( \ph{\h{Z}_{n}}{t\wedge \tstar{\h{Z}_{n}} } ) 
\1_{\{\h{S}_{n+1} \geq t\wedge \tstar{\h{Z}_{n}}\}} \big| \h{Z}_{n}= z\Big], \\
\h{L}^{d}_{n+1}(v,w)(z) & = \min_{t \in G(z)} \h{J}_{n+1}(v,w)(z,t) \wedge \h{K}_{n+1}w(z).
\end{align*}
\end{definition}
With these discretized operators and the previous evaluation of the $\widetilde{v}_k$, we propose the following approximation scheme.
\begin{definition}
\label{vchap}
Consider $\widehat{v}_{N}(z)=g(z)$ where $z\in \Gamma^{Z}_{N}$ and for $k\in\{1,\ldots,N\}$
\begin{eqnarray}
\widehat{v}_{k-1}(z) & = & \widehat{L}_{k}^{d}(M\widetilde{v}_{k},\widehat{v}_{k})(z), 
\end{eqnarray}
where $z\in \Gamma^{Z}_{k-1}$.
\end{definition}
See figure figure \ref{algo} for a graphical illustration of this numerical procedure.
Therefore $\widehat{v}_0(\widehat{Z}_0)$ will be an approximation of $v_0(Z_0)=v_0(x_0)$. The derivation of the error bound for this scheme follows exactely the same lines as in the preceding section. Therefore we omit it and only state our main result.

\begin{theorem}
\label{th2 erreur control}
For all $0\leq n\leq N-1$, suppose that $\Delta(z)$ for $z\in \Gamma^{Z}_{n}$ is such that
\begin{eqnarray*}
\sqrt{\frac{D^{4} \|Z_{n}-\h{Z}_{n}\|_{p} + D^{5} \big\| S_{n+1}- \h{S}_{n+1} \big\|_{p}}{D^{3}}}&<& \min_{z\in \Gamma^{Z}_{n}}\{\Delta(z)\}.
\end{eqnarray*}
Then we have
\begin{eqnarray*}
\lefteqn{\|v_{n}(Z_n)-\widehat{v}_{n}(\h{Z}_n)\|_p}\\
& \leq & \big\|v_{n+1}(Z_{n+1})-\h{v}_{n+1}(\widehat{Z}_{n+1})\big\|_p
+ \max_{y\in \mathbb{U}} \big| v_{n+1}(y)-\widetilde{v}_{n+1}(y)\big|
+  D^{1}_{n} \|Z_{n}-\widehat{Z}_{n}\|_p\\
&& +3\Li{v_{n+1}} \big\| Z_{n+1}- \h{Z}_{n+1} \big\|_{p}
+2C_{f} \big\| S_{n+1}- \h{S}_{n+1} \big\|_{p} +D^{2}_{n} \big\| \Delta(\widehat{Z}_{n}) \big\|_p \\
&& +2 \sqrt{D^{3} \big( D^{4} \|Z_{n}-\h{Z}_{n}\|_{p} + D^{5} \big\| S_{n+1}- \h{S}_{n+1} \big\|_{p} \big)},
\end{eqnarray*}
\end{theorem}
with
\begin{eqnarray*}
D^{1}_{n} & = & [v_{n}]+ \Li{Q} \Li{v_{n+1}}_1\frac{C_{\lambda}}{\alpha} +\frac{C_{f}}{\alpha}\big(E_{1}+E_{2}\big) \\
&& + \Big\{ \Li{Q} \Li{v_{n+1}}_{*}  + 2 E_{3} \Big\} \vee \Big\{ 2 \big( \Li{c}_{1}+\Li{c}_{2}\Li{t^{*}} \big)
+C_{c} E_{1} + \alpha \Li{t^{*}} (\frac{C_{f}}{\alpha}+C_{c}) 
\Big\} ,\\
D^{2}_{n} & = &  C_f + C_{v_{n+1}}C_{\lambda} + \Li{c}_2+(C_c+C_{v_{n+1}})(C_{\lambda}+\alpha ), \\
D^{3} & = & \big(\frac{2C_{f}}{\alpha}+C_{c} \big) C_{\lambda},\\
D^{4} & = & 2 \Li{t^{*}} \big( \frac{2C_{f}}{\alpha} +C_{c} \big), \\
D^{5} & = & 2 \big( \frac{2C_{f}}{\alpha} +C_{c} \big).
\end{eqnarray*}
\begin{remark}
By using the same arguments as in Remark \ref{larochelle}, it can be shown that
$\|v_{n}(Z_n)-\widehat{v}_{n}(\h{Z}_n)\|_p$
can be made arbitrarily small by an adequate choice of the discretization parameters.
From a theoretical point of view, the error can be calculated by iterating the result of Theorem \ref{th2 erreur control}.
However, this result is not presented here because it would lead to an intricate expression. From a numerical point of view, a computer  can easily estimate this error as shown in the example of section \ref{section appli}.
\end{remark}

\subsection{Step by step description of the algorithm}
\begin{figure}[htbp]
\begin{center}
\fbox{\includegraphics[width=23cm,angle=-90]{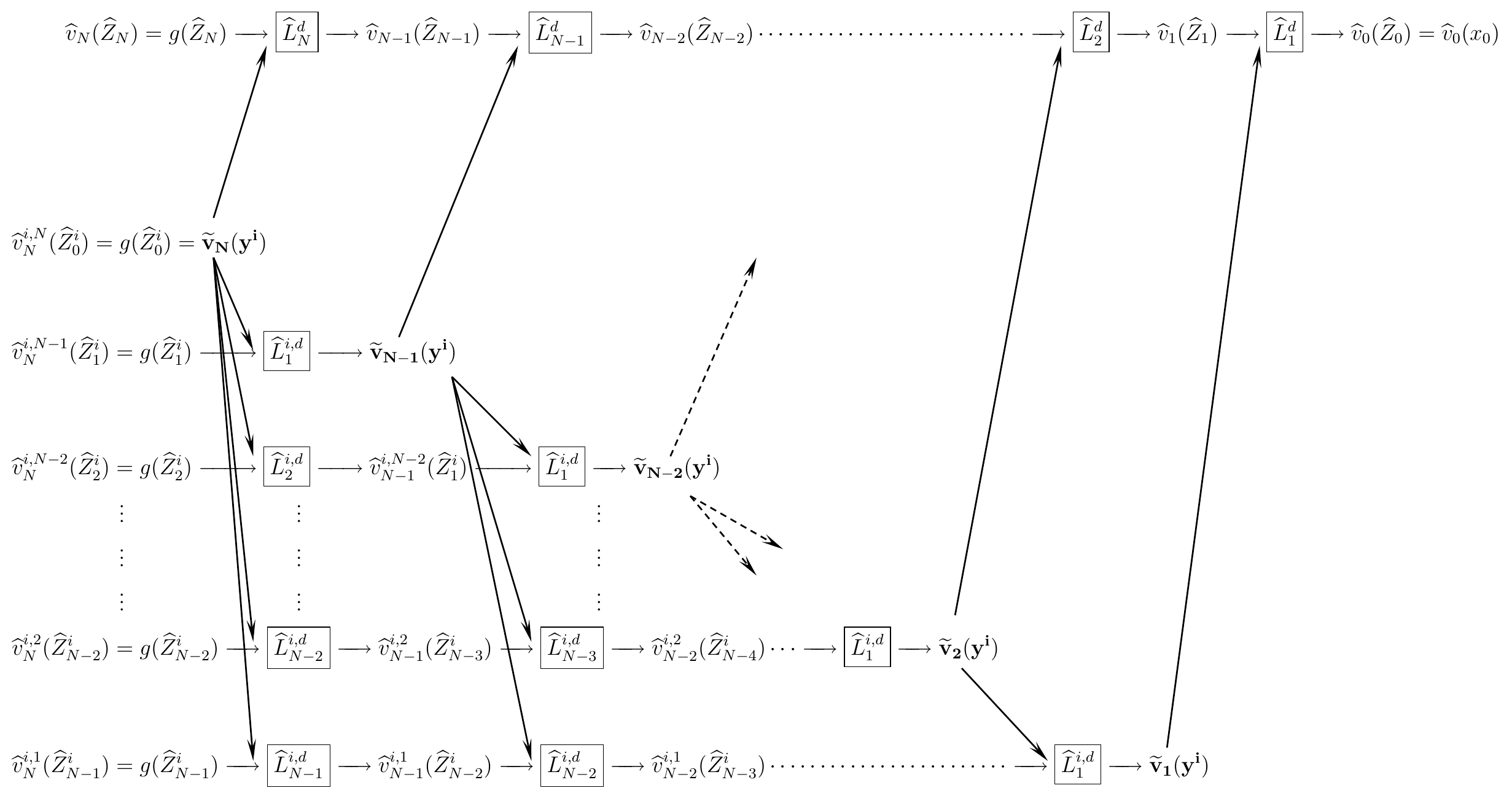}}
\caption{Step by step procedure}
\label{algo}
\end{center}
\end{figure}

Recall that the main objective of our algorithm is to compute the approximation $\widehat{v}_{0}(x_{0})$ of the value function of the impulse control problem $v_{0}(x_{0})$. 
The global recursive procedure is described on figure \ref{algo}. 

\bigskip

\noindent
The calculation of $\widehat{v}_{0}(x_{0})$ is based on the backward recursion given in Definition \ref{vchap} and described in the first line of figure \ref{algo}. It involves the operators $\widehat{L}^{d}_{j}$ constructed with the quantized process $\widehat{\Theta}_{n}$ starting from $x_{0}$. This recursion is not self contained and requires previous evaluation of the functions $\widetilde{v}_{j}$ on the control set $\mathbb{U}$.

\bigskip

\noindent
The lower part of figure \ref{algo}, shows how to compute these functions $\widetilde{v}_{j}$ at each point of the control grid
$\mathbb{U}$. This is the triangular backward recursion given in Definition \ref{vchapi}. 
More precisely, define $\widetilde{v}_{N}=g$ and set $j < N$ and suppose that  all the $\widetilde{v}_{l}$ for all $l>j$ have already been computed everywhere on the control set $\mathbb{U}$. One then computes $\widetilde{v}_{j}$ in the following way, following the $j$-th line of figure  \ref{algo} counting from the bottom. One first iterates the operators 
$\widehat{L}^{1,d}_{k}$ and uses the quantized process $\widehat{\Theta}^{1}_{n}$, to obtain $\widetilde{v}_{j}(y^{1})$. Then one iterates the operators 
$\widehat{L}^{2,d}_{k}$ and uses the quantized process $\widehat{\Theta}^{2}_{n}$, to obtain $\widetilde{v}_{j}(y^{2})$, and so on until the last point $\widetilde{v}_{j}(y^{u})$.
Thus one obtains $\widetilde{v}_{j}$ at all points of the control set $\mathbb{U}$.

\subsection{Practical implementation}
The procedure defined above is the natural one to obtain convergence rates for our approximations. However, in practice we proceed in a different order.

\bigskip

The first step is to fix the computational horizon $N$. This point was discussed earlier. The second step is not the time discretization, but the computation of the quantized approximations of the sequences $(\Theta_n)$ and $(\Theta_n^i)$. The quantization algorithm may be quite long to run. However, it must be pointed out that this quantization step only depends on the optimization procedure through the
the control set $\mathbb{U}$ but  it does not depend on the cost functions $f$ and $c$.
The sequence $(\widehat{\Theta}_n)$ is obtained in a straightforward way. As for the $(\widehat{\Theta}_n^i)$, if the control set is very small, it is possible to run as many sequences of grids as there are points in the control set. Otherwise, one can do with only one sequence of grids computed with the Markov chain $(\Theta_n^{\mu})$ with a random starting point $Z_0=Z_0^{\mu}$ uniformly distributed on the control set $\mathbb{U}$. To derive the point-wise approximation error, one simply uses the finiteness of $\mathbb{U}$ and the definition of the $L_p$ norm.
\begin{eqnarray*}
|v_k(y^i)-\widetilde{v}_k(y^i)|&\leq& u\sum_{i=1}^u|v_k(y^i)-\widetilde{v}_k(y^i)|\frac{1}{u}\\
&\leq&u^{1/p}\|v_k(Z_0^{\mu})-\widetilde{v}_k(\widehat{Z}_0^{\mu})\|_p
\end{eqnarray*}
where $u$ is the cardinal of $\mathbb{U}$. Notice that the last term is bounded in Theorem~\ref{th erreur control}. Hence, one really only needs two series of quantization grids.

\bigskip

Once the quantization grids are computed and stored, one computes the path-adapted time grids $G(z)$ for all $z$ in all the quantization grids, that is only a finite number of $z$. The step $\Delta(z)$ can usually be chosen constant equal to $\Delta$, so that either one can store the whole time grids, or one only needs to store the values of $\Delta$ and $t^*(z)$ for all $z$ in the quantization grids.

\bigskip

Once these preliminary computations are done, one can finally compute the value function. This last step is comparatively faster. The only point left to discussion is how to choose the initializing function $g$. The most interesting starting function is the cost $h$ of the no impulse strategy, because then the value function $\mathcal{V}_N$ has a natural interpretation. However, in general, one needs additional assumptions on $Q$ to ensure that $h$ is in $\mathbf{L}$. Another problem, is that in general computing $h$ is a difficult problem, especially as we need to know its value at many different points, as explained in Remark~\ref{remark g}. To overcome these difficulties, one can choose $g$ to be an upper bound of $h$, for instance, $g=\alpha^{-1}C_f$. In the special cases where $h$ can be explicitly computed, we advise to use $h$.

\section{Example}
\label{section appli}
Now we apply our procedure to a simple PDMP and present numerical results. This example is quite similar to example (54.29) in \cite{davis93}, we only added random jumps to obtain a non trivial Markov chain $(Z_n,S_n)$.

\bigskip

\begin{figure}[htbp]
\begin{center}
\begin{minipage}{.45\textwidth}
\begin{center}
\includegraphics[width=7.25cm]{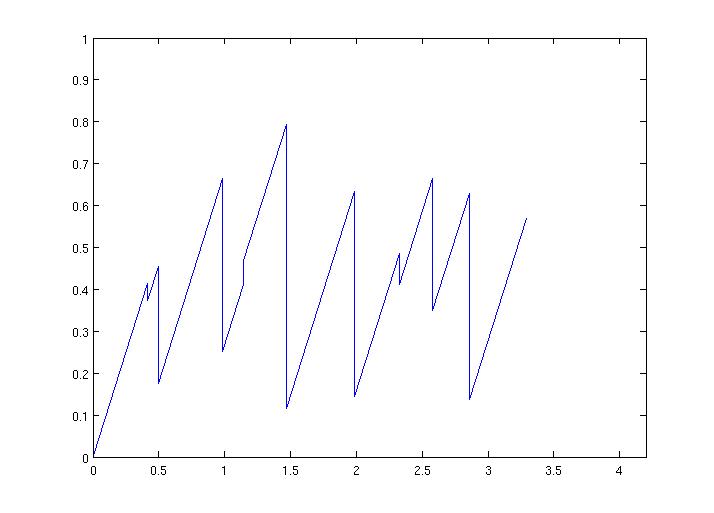}
\end{center}
\end{minipage}
\begin{minipage}{.45\textwidth}
\begin{center}
\includegraphics[width=7cm]{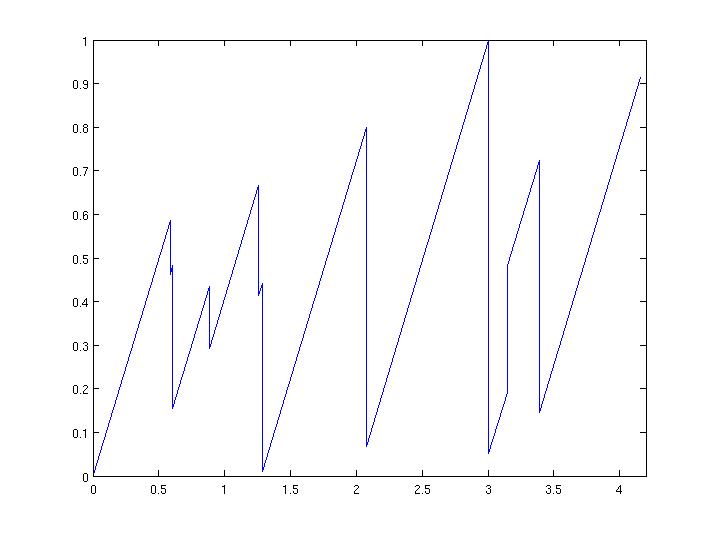}
\end{center}
\end{minipage}
\caption{Two trajectories of the PDMP.}\label{figure traj}
\end{center}
\end{figure}

Set $E=[0, 1[$, and $\partial E=\{1\}$. The flow is defined on $[0,1]$ by $\phi(x,t)=x+vt$ for some positive $v$, the jump rate is defined on $[0,1]$ by $\lambda(x)=\beta x$, with $\beta>0$, and for all $x\in [0,1]$, one sets $Q(x,\cdot)$ to be the uniform law on $[0, 1/2]$. Thus, the process moves with constant speed $v$ towards $1$, but the closer it gets to the boundary $1$, the higher the probability to jump backwards on $[0, 1/2]$. Figure~\ref{figure traj} shows two trajectories of this process for $x_0=0$, $v=1$ and $\beta=3$ and up to the $10$-th jump.
The running cost is defined on $E$ by $f(x)=1-x$ and the intervention cost is a constant $c_0$. Therefore, the best performance is obtained when the process is close to the boundary $1$. The control set $\mathbb{U}$ is the set of $\frac{k}{u}$, $0\leq k\leq u-1$ for some fixed integer $u$. In this special case, the control grid is already a discretization of the whole state space of the process. Therefore one needs only one series of grids starting from the control points to obtain an approximation of the value function at each point of the control grid. 
\begin{figure}[p]
\begin{center}
\includegraphics[width=8cm]{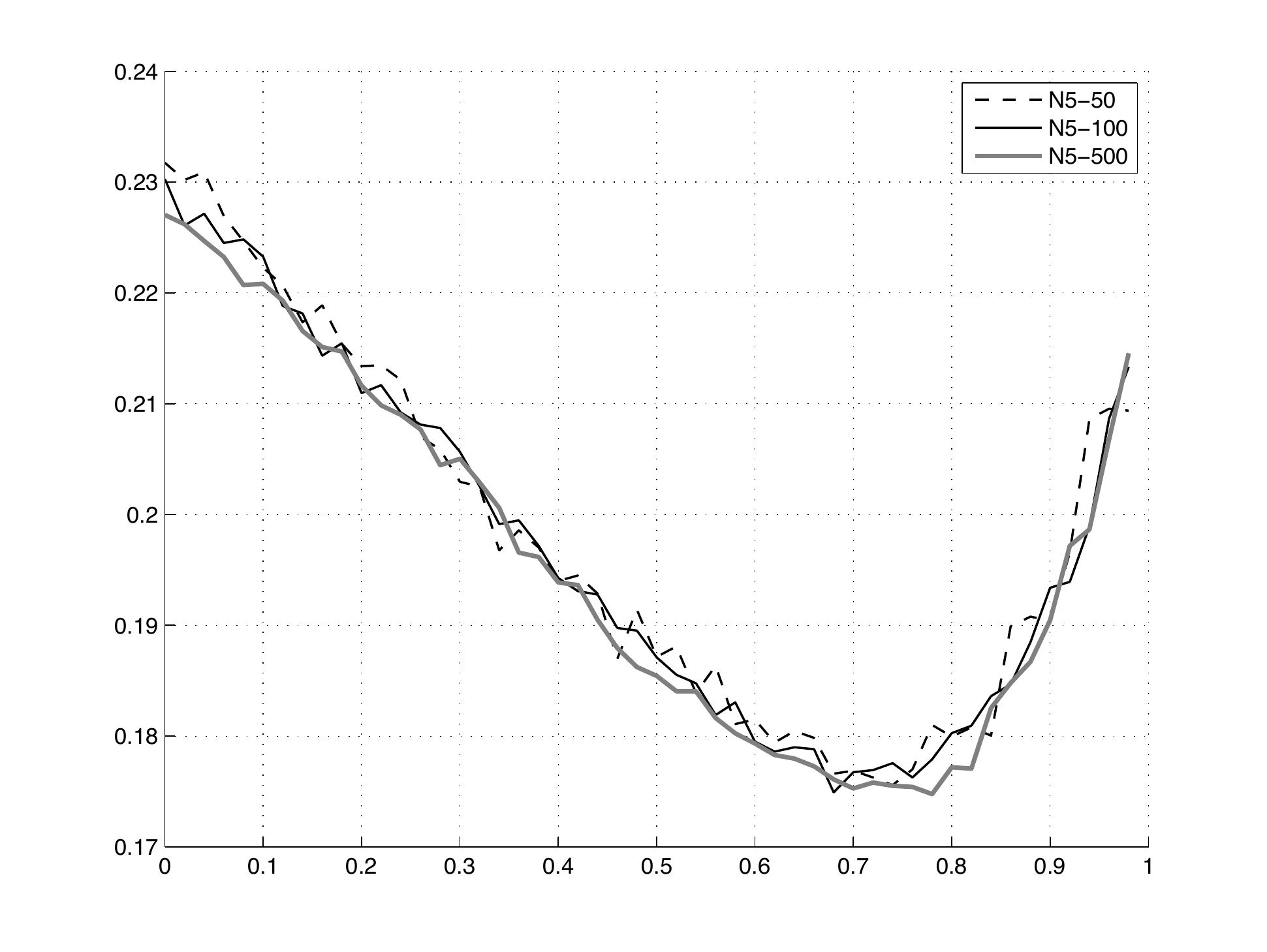}
\caption{Approximated value function for $N=5$.}\label{figure value 5}
\end{center}
\end{figure}

\begin{figure}[p]
\begin{center}
\includegraphics[width=8cm]{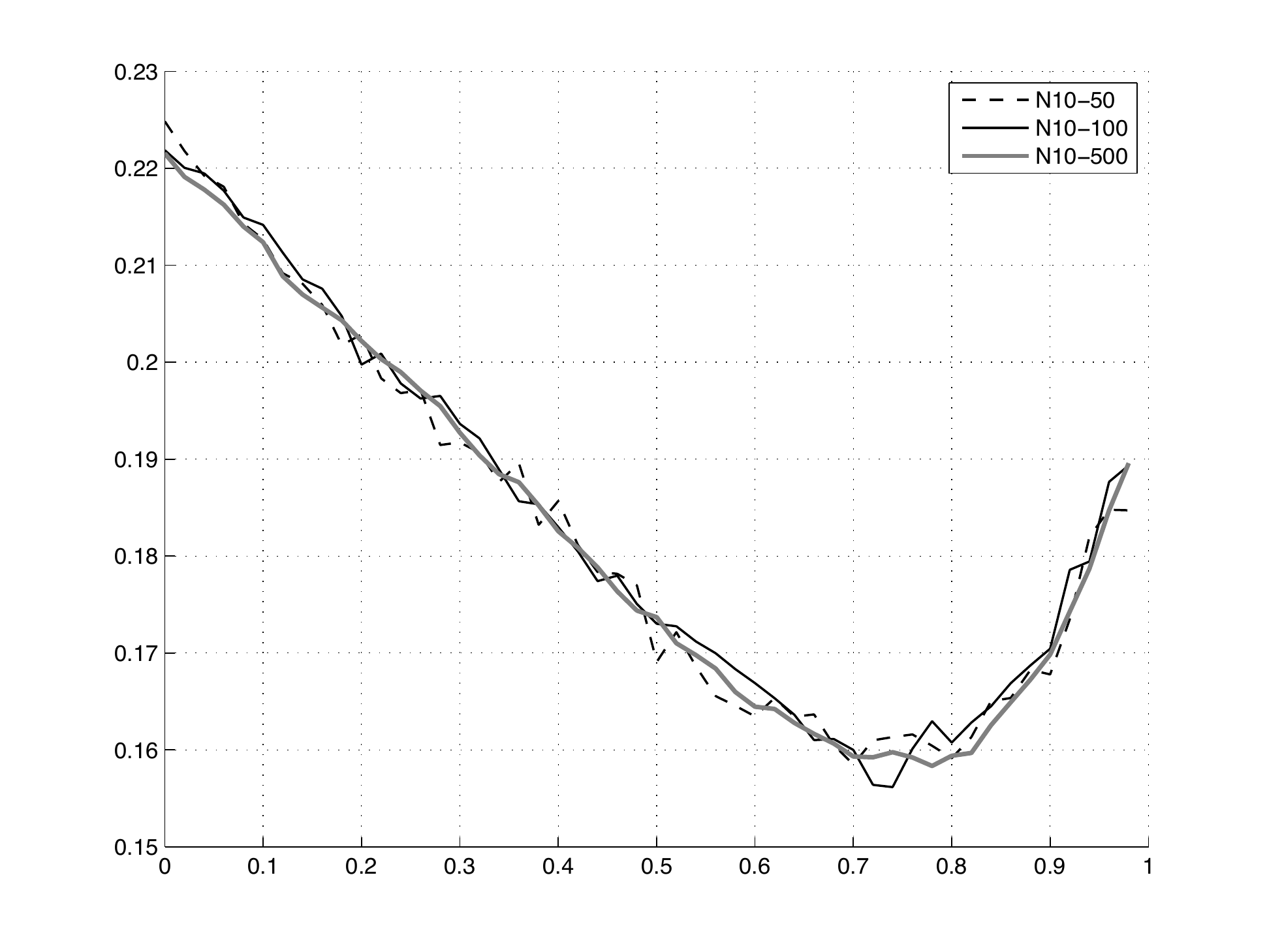}
\caption{Approximated value function for $N=10$.}\label{figure value 10}
\end{center}
\end{figure}

\begin{figure}[p]
\begin{center}
\includegraphics[width=8cm]{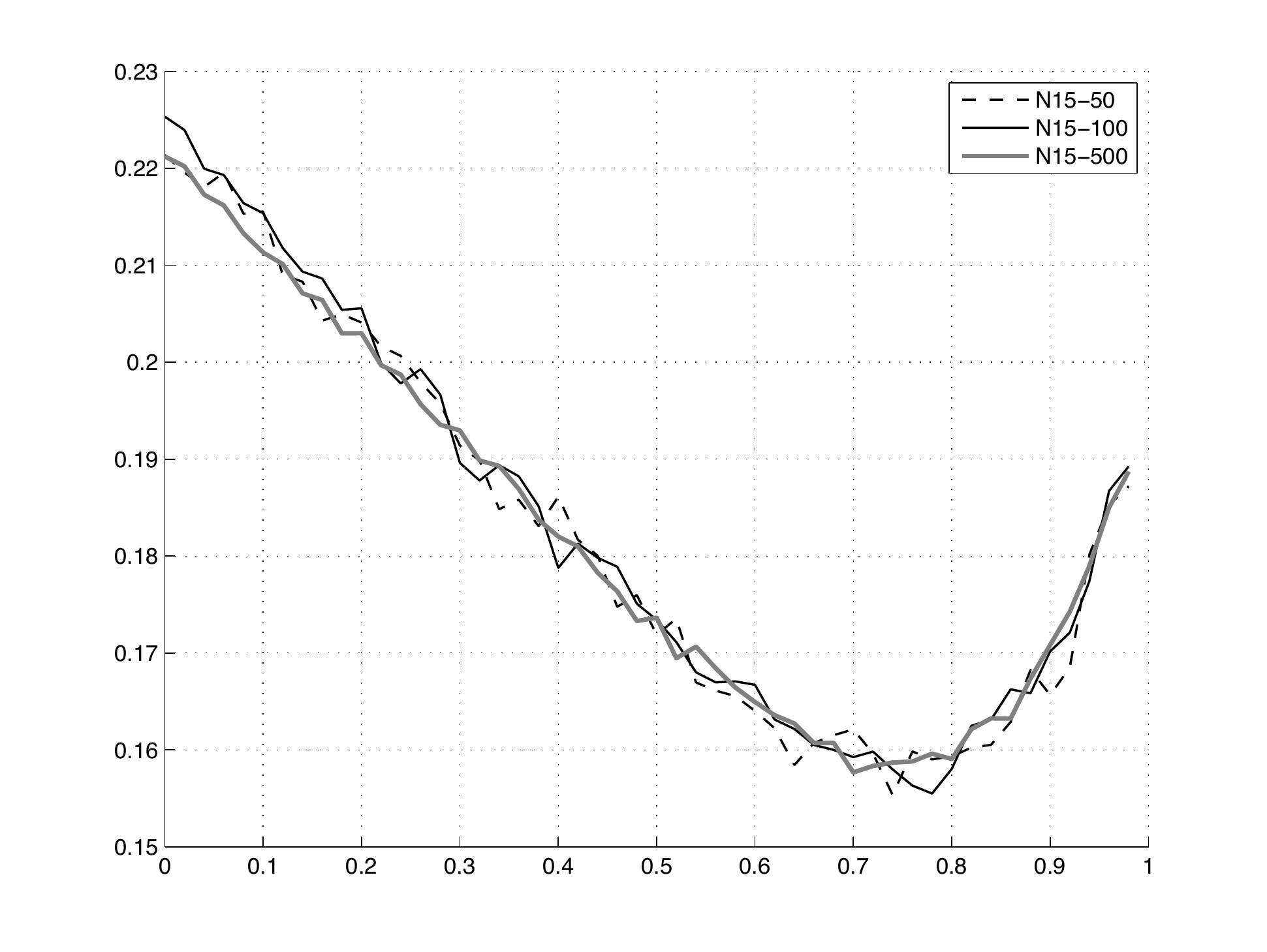}
\caption{Approximated value function for $N=15$.}\label{figure value 15}
\end{center}
\end{figure}
We ran our algorithm for the parameters  $x_0=0$, $v=1$, $\beta=3$, $c_0=0.08$, the discount factor $\alpha=2$ and $u=50$ points in the control grid and several values of the horizon $N$. 

\bigskip

For an horizon $N=5$ (respectively, $N=10$, $N=15$) interventions or jumps, Figure~\ref{figure value 5} (respectively, Figure~\ref{figure value 10}, Figure~\ref{figure value 15})
gives the approximated value function we obtained (computed at the 50 points of the 
control grid) for $50$, $100$ and $500$ discretization points in each quantization grid and. As expected, the approximation gets smoother and lower as the number of points in the quantization grids increases. 

\bigskip

The theoretical errors corresponding to the horizon $N=5$ (respectively, $N=10$, $N=15$) are given in Table~\ref{table5} (respectively, Table~\ref{table10}, Table~\ref{table15}).
The values of the error are fairly high and conservative, but it must be pointed out that on the one hand, they do decrease as the number of points in the quantization grids increase, as expected ; on the other hand their expressions are calculated and valid for a very wide and general class of PDMP's, hence when applied to a specific example, they cannot be very sharp.
\begin{table}[hp]
\begin{center}
\begin{tabular}{|c|c|}
\hline
Number of points in the quantization grids&$\|v_0(Z_0)-\widehat{v}_0(Z_0)\|_2$\\
\hline
$50$&4636\\
$100$&3700\\
$500$&2141\\
\hline
\end{tabular}
\caption{Theoretical errors for $N=5$.}
\label{table5}
\end{center}
\end{table}

\begin{table}[hp]
\begin{center}
\begin{tabular}{|c|c|}
\hline
Number of points in the quantization grids&$\|v_0(Z_0)-\widehat{v}_0(Z_0)\|_2$\\
\hline
$50$&5.341$\cdot10^{11}$\\
$100$&4.501$\cdot10^{11}$\\
$500$&2.567$\cdot10^{11}$\\
\hline
\end{tabular}
\caption{Theoretical errors for $N=10$.}
\label{table10}
\end{center}
\end{table}

\begin{table}[hp]
\begin{center}
\begin{tabular}{|c|c|}
\hline
Number of points in the quantization grids&$\|v_0(Z_0)-\widehat{v}_0(Z_0)\|_2$\\
\hline
$50$&1.460$\cdot10^{21}$\\
$100$&1.288$\cdot10^{21}$\\
$500$&0.750$\cdot10^{21}$\\
\hline
\end{tabular}
\caption{Theoretical errors for $N=15$.}
\label{table15}
\end{center}
\end{table}


\bigskip

\begin{figure}[htp]
\begin{center}
\includegraphics[width=10cm]{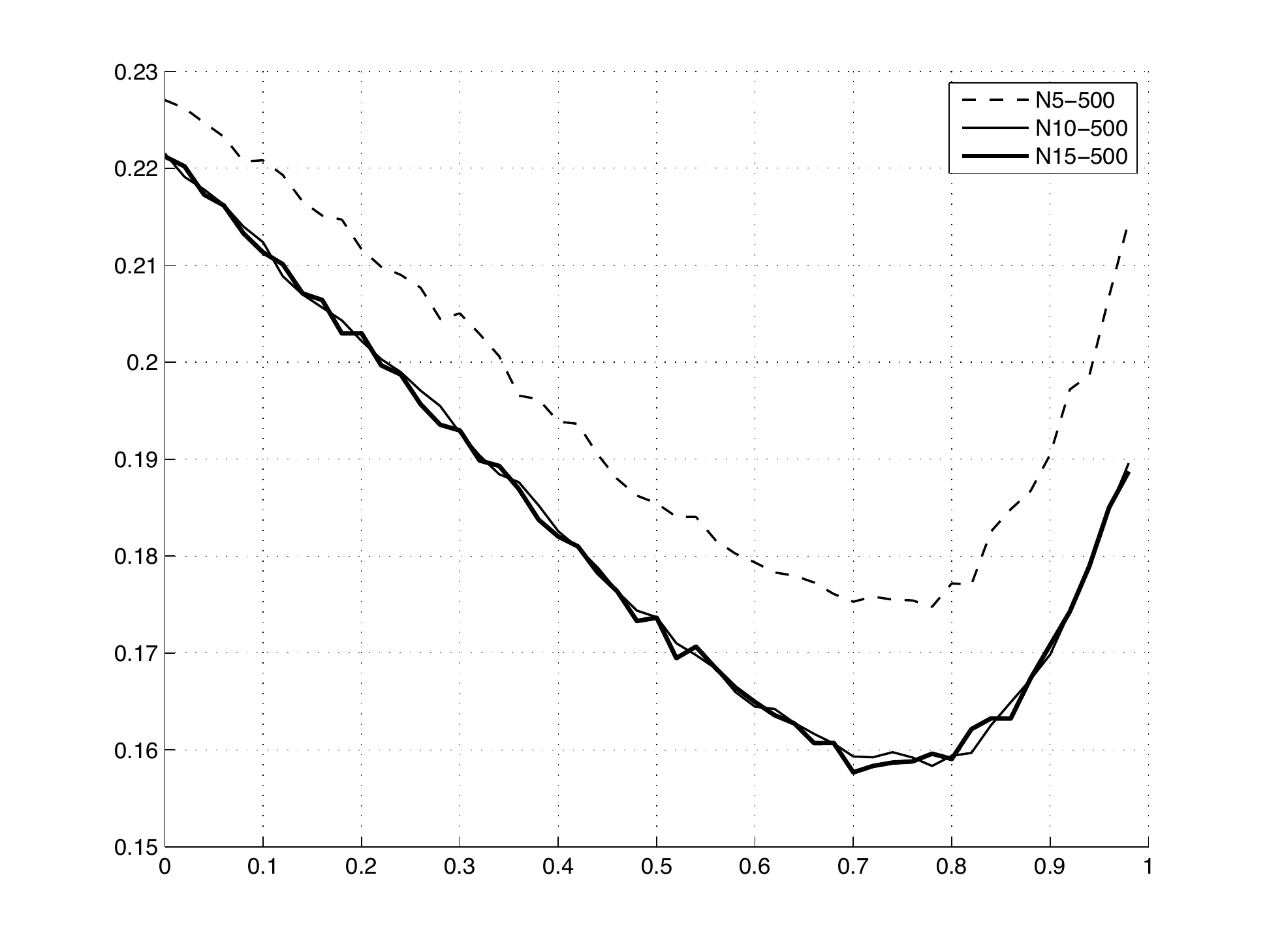}
\caption{Approximated value function for $N=5$, $N=10$ and $N=15$ for $500$ points in the quantization grids.}\label{figure compare500}
\end{center}
\end{figure}

Notice also that the approximated value function obtained for the horizon of $N=10$ jumps or interventions is much lower than that obtained for the horizon $N=5$ jumps or interventions. This is natural as it is a minimization problem, and the more there are possible interventions the lower the value function is. This also suggests that the horizon $N=5$ is not large enough to approximate the infinite horizon problem.   
Figure~\ref{figure compare500} gives the approximated value function we obtained (computed at the 50 points of the 
control grid) for respectively $500$ points in the quantization grids and respective horizons of $N=5$, $N=10$ and $N=15$ jumps or interventions. There is very little difference between the results for $N=10$ and $N=15$, suggesting that it is enough to take an horizon of $10$ jumps or intervention to approximate the infinite time optimization problem.

\appendix
\section{Lipschitz continuity results}
\subsection{Lipschitz properties of the operators}
We start with preliminary results on operators $M$, $H$, $F$ and $I$.
\begin{lemma}\label{lem M}
For any $g$ defined on $\mathbb{U}$, $Mg\in \mathbf{L}$. Moreover,
$$\Li{Mg}_{1}\leq \Li{c}_{1}, \qquad \Li{Mg}_{2} \leq \Li{c}_{2}, \qquad \Li{Mg}_{*} \leq \Li{c}_{*},\qquad C_{Mg}\leq C_c+C_g.$$
\end{lemma}

\noindent{\textbf{Proof:}} 
By using the fact that $\big|Mg(x)-Mg(y)\big| \leq \sup_{z\in \mathbb{U}} \big| c(x,z)-c(y,z)\big|$ and assumption~\ref{H7}, the result follows easily.
\hfill $\Box$

\begin{lemma}
\label{TechH}
Let $v\in \mathbf{L}$. Then for all $(x,y)\in E^{2}$ and $(t,u)\in \mathbb{R}_+^2$, one has
\begin{equation*}
\Big|Hv(x,t)-Hv(y,u)|\leq D_{1}(v) |x-y| +D_{2}(v) |t-u|,
\end{equation*}
where
\begin{itemize}
\item if $t<\tstar{x}$ and $u<\tstar{y}$,
\begin{equation*}
D_1(v)=\Li{v}_{1}+C_{v}C_{t^{*}}\Li{\lambda}_{1},\qquad D_2(v)=\Li{v}_2+C_{v}( C_{\lambda}+\alpha),
\end{equation*}
\item if $t=\tstar{x}$ and $u=\tstar{y}$,
\begin{equation*}
D_1(v)=\Li{v}_*+C_{v}\big(C_{t^{*}}\Li{\lambda}_1+(C_{\lambda}+\alpha)\Li{t^{*}}\big),\qquad D_{2}(v)=0,
\end{equation*}
\item otherwise,
\begin{equation*}
D_1(v)=\Li{v}_{1}+\Li{v}_{2}\Li{t^{*}}+C_{v}\big(C_{t^{*}}\Li{\lambda}_1+(C_{\lambda}+\alpha) \Li{t^{*}}\big),\qquad D_2(v)=\Li{v}_{2}+C_{v} (C_{\lambda}+\alpha).
\end{equation*}
\end{itemize}
\end{lemma}

\noindent{\textbf{Proof:}} 
This is straightforward.
\hfill $\Box$

\begin{lemma}
\label{TechFI}
For all $w\in \mathbf{L}$, $(x,y)\in E^{2}$ and $(t,u)\in \mathbb{R}_+^2$, one has
\begin{eqnarray*}
\Big|F(x,t)-F(y,u)\Big|
&\leq&\frac{1}{\alpha}\Big( \Li{f}_1+C_{f} C_{t^{*}}\Li{\lambda}_1 \Big)|x-y|+C_{f}\big(|t-u|\vee\Li{t^{*}}|x-y|\big),\\
\Big|Iw(x,t)-Iw(y,u)\Big| &\leq&
\frac{1}{\alpha} \Big( \Li{Q}  \Li{w}_1 C_{\lambda} +C_w\Li{\lambda}_1 \big(1+C_{\lambda}C_{t^{*}}\big) \Big) |x-y|\\
&&+C_wC_{\lambda}\Big(|t-u|\vee\Li{t^{*}}|x-y|\Big).
\end{eqnarray*}
\end{lemma}

\noindent{\textbf{Proof:}} 
Suppose, without loss of generality, that $t\wedge\tstar{x}\leq u\wedge\tstar{y}$. Then, one has
\begin{align*}
\Big|F(x,t)-F(y,u)\Big|
\leq&\int_{0}^{t\wedge\tstar{x}}\mathrm{e}^{-\alpha s}\Big|f\big(\ph{x}{s} \big)\expl{x}{s} - f\big(\ph{y}{s}\big)\expl{y}{s}\Big|ds\\
&+\int_{t\wedge\tstar{x}}^{u\wedge\tstar{y}}\Big|f\big(\ph{y}{s} \big)\expla{y}{s}\Big|ds\\
\leq&\Big( \Li{f}_{1}+C_{f}C_{t^{*}}\Li{\lambda}_{1} \Big) \int_0^{\infty}\mathrm{e}^{-\alpha s}ds |x-y|+C_f\big|u\wedge\tstar{y}-t\wedge\tstar{x}\big|.
\end{align*}
From the fact that $\big|u\wedge\tstar{y}-t\wedge\tstar{x}\big|\leq |t-u|\vee\Li{t^{*}}|x-y|$ we get the first inequality.
\nl
By using similar arguments, it is easy to obtain the last result.
\hfill $\Box$

Now we turn to the Lipschitz property of operator $K$.

\begin{lemma}\label{TechK}
For $w\in \mathbf{L}$ and $(x,y)\in E^2$, one has
\begin{align*}
\Big|Kw(x)-Kw(y) \Big| \leq& \bigg\{
\Li{Q} \Li{w}_1\frac{C_{\lambda}}{\alpha}+\Li{Q} \Li{w}_{*} +C_w\Big(E_{1}+E_{2}\Big) +E_{3} \bigg\}|x-y|.
\end{align*}
\end{lemma}

\noindent{\textbf{Proof:}} 
This is a direct consequence of (\ref{defK}) and Lemmas~\ref{TechH}, \ref{TechFI}.
\hfill $\Box$

Finally, we study the the Lipschitz properties of operator $J$.

\begin{lemma}
\label{TechJt}
For all $(v,w)\in \mathbf{C}^{2}$, $x\in E$ and $(t,u)\in \mathbb{R}_+^2$, one has
\begin{equation*}
\Big|J(v,w)(x,t)-J(v,w)(x,u)\Big|\leq\Big( C_f + C_wC_{\lambda} + \Li{v}_2+C_v(C_{\lambda}+\alpha)\Big)|t-u|.
\end{equation*}
\end{lemma}

\noindent{\textbf{Proof:}} 
By using (\ref{defJ}) and Lemmas~\ref{TechH} and \ref{TechFI}, the result follows easily.
\hfill $\Box$

\begin{lemma}
\label{TechJx}
For all $(v,w)\in \mathbf{L}^{2}$, $(x,y)\in E^2$ and $t\geq0$, one has
\begin{align*}
\Big|J(v,w)(x,t)-J(v,w)(y,t)\Big|  \leq 
& \Big\{ \Li{v}_1+\Li{v}_2\Li{t^{*}}+C_v  E_{1}+\Li{Q}\Li{w}_1\frac{C_{\lambda}}{\alpha} \\
& +C_w E_{2}+E_{3} \Big\} |x-y|.
\end{align*}
where
\begin{eqnarray*}
E_{1} & = & C_{t^{*}}\Li{\lambda}_1+(C_{\lambda}+\alpha) \Li{t^{*}}, \\
E_{2} & = & C_{\lambda}\Li{t^{*}}+\Li{\lambda}_1\frac{1+C_{\lambda}C_{t^{*}}}{\alpha}, \\
E_{3} & = & \Li{f}_1\frac{1}{\alpha}+C_f \Big(\frac{C_{t^{*}}\Li{\lambda}_1}{\alpha}+\Li{t^{*}}\Big).
\end{eqnarray*}
\end{lemma}

\noindent{\textbf{Proof:}} 
Again, this is a direct consequence of (\ref{defJ}) and Lemmas~\ref{TechH} and \ref{TechFI}.
\hfill $\Box$

\begin{remark}
\label{InfJx}
It is easy to obtain
that for $(v,w)\in \mathbf{C}^{2}$, $s \in \RR_{+}$ and $(x,y) \in E^{2}$,
\begin{align*}
\Big| \inf_{t\geq s}J(v,w)(x,t) - \inf_{t\geq s}J( v,w)(y,t) \Big| & \leq 
\sup_{t\geq 0} \big| J(v,w)(x,t) -  J(v,w)(y,t) \big|
\end{align*}
\end{remark}

\begin{lemma}
\label{TechInfJt}
Let $(v,w)\in \mathbf{L}^{2}$. Then for all $x \in E$ and $(s,t)\in \RR_{+}^{2}$,
\begin{eqnarray*}
\Big| \inf_{u\geq t}J(v,w)(x,u) - \inf_{u\geq s}J(v,w)(x,u) \Big| \leq \Big( C_f + C_wC_{\lambda} + \Li{v}_2+C_v(C_{\lambda}+\alpha)\Big)
|t-s| .
\end{eqnarray*}
\end{lemma}

\noindent{\textbf{Proof:}} 
Without loss of generality it can be assumed that $s\leq t$.
Therefore, one has
\begin{align}
\Big| \inf_{u\geq t}J(v,w)(x,u) - \inf_{u\geq s}J(v,w)(x,u) \Big| =
\inf_{u\geq t}J(v,w)(x,u) - \inf_{u\geq s}J(v,w)(x,u).
\label{eqT1}
\end{align}
Remark that there exists $\widebar{s} \in [s\wedge \tstar{x},\tstar{x}]$ such that
$\ds \inf_{u\geq s}J(w,g)(x,u) = J(w,g)(x,\widebar{s})$.
Consequently, if $\widebar{s}\geq t\wedge \tstar{x}$ then one has $\ds \Big| \inf_{u\geq t}J(v,w)(x,u) - \inf_{u\geq s}J(v,w)(x,u) \Big|=0$.\\
Now if $\widebar{s} \in [s\wedge \tstar{x},t\wedge \tstar{x}[$, then one has
\begin{eqnarray*}
\inf_{u\geq t}J(v,w)(x,u) - \inf_{u\geq s}J(v,w)(x,u) \leq J(v,w)(x,t) - J(v,w)(x,\widebar{s}).
\end{eqnarray*}
From Lemma \ref{TechJt}, we obtain the following inequality
\begin{eqnarray}
\inf_{u\geq t}J(v,w)(x,u) - \inf_{u\geq s}J(v,w)(x,u) \leq
 \Big( C_f + C_wC_{\lambda} + \Li{v}_2+C_v(C_{\lambda}+\alpha)\Big) |t-\widebar{s}|.
\label{eqT2}
\end{eqnarray}
Combining equations (\ref{eqT1}), (\ref{eqT2}) and the fact that  $|t-\widebar{s}| \leq |t-s|$ the result follows.
\hfill $\Box$

\subsection{Lipschitz properties of the operator $\mathcal{L}$}
Now we study the Lipschitz continnuity of our main operator
\begin{lemma}
\label{TechFIHphi}
For all $(v,w)\in \mathbf{L}^{2}$, $x\in E$ and $t\in [0,\tstar{x})$ and $u\in \mathbb{R}_{+}$, one has
\begin{eqnarray*}
F\big(\ph{x}{t}, u\big) & = & \mathrm{e}^{\alpha t+\Lambda(x,t)}\Big(F(x,t+u)-F(x,u)\Big), \\
Iw\big(\ph{x}{t}, u\big) & = & \mathrm{e}^{\alpha t+\Lambda(x,t)}\Big(Iw(x,t+u)-Iw(x,u)\Big), \\
Hv\big(\ph{x}{t}, u\big) & = & \mathrm{e}^{\alpha t+\Lambda(x,t)}Hv(x,t+u).
\end{eqnarray*}
\end{lemma}

\noindent{\textbf{Proof:}} 
By using the semi-group property of $\phi$, we have $\Lambda\big(\ph{x}{t}, u\big)=\Lambda(x,t+u)-\Lambda(x,t)$ for
$t+u<\tstar{x}$ and noting that $\tstar{\ph{x}{t}}=\tstar{x}-t$  for $t<\tstar{x}$, a simple change of variable yields
\begin{equation*}
F\big(\ph{x}{t}, u\big)=\mathrm{e}^{\alpha t+\Lambda(x,t)}\int_t^{(t+u)\wedge \tstar{x}}\expla{x}{s} f \big(\ph{x}{s}\big)ds,
\end{equation*}
and we get the first equation.
The other equalities can be obtained by using similar arguments.
\hfill $\Box$

\begin{lemma}
\label{TechLphi}
For all $(v,w)\in \mathbf{L}^{2}$, $x\in E$ and $t\in [0,\tstar{x})$,
\begin{eqnarray*}
L(v,w)\big(\ph{x}{t}\big) & = & \mathrm{e}^{\alpha t+\Lambda(x,t)}
\bigg[
\Big\{ \inf_{s\geq t} J(v,w) (x,s) \wedge Kw(x) \Big\} -F(x,t)- Iw(x,t)
\bigg].
\end{eqnarray*}
\end{lemma}

\noindent{\textbf{Proof:}} 
For $t\in [0,\tstar{x})$ and $u\in \mathbb{R}_{+}$, we have from Lemma \ref{TechFIHphi}, (\ref{defJ}) and (\ref{defK})
\begin{eqnarray*}
J(v,w)\big(\ph{x}{t},u\big)&=&\mathrm{e}^{\alpha t+\Lambda(x,t)}\big[J(v,w)(x,t+u)-F(x,t)-Iw(x,t)\big],\\
Kw\big(\ph{x}{t}\big)&=&\mathrm{e}^{\alpha t+\Lambda(x,t)}\big[Kw(x)-F(x,t)-Iw(x,t)\big].
\end{eqnarray*}
Consequently, from equation (\ref{DefL}), it follows
\begin{eqnarray*}
L(v,w)\big(\ph{x}{t}\big) & = & \mathrm{e}^{\alpha t+\Lambda(x,t)}
\bigg[
\Big\{ \inf_{u\geq 0} J(v,w) (x,t+u) \wedge Kw(x) \Big\} -F(x,t)- Iw(x,t)
\bigg],
\end{eqnarray*}
showing the result.
\hfill $\Box$

\begin{proposition}
\label{TechLiL}
For all $w\in \mathbf{L}$, $\mathcal{L}w\in  \mathbf{L}$, and
\begin{eqnarray*}
\Li{\mathcal{L}w}_{1} & \leq & \mathrm{e}^{(\alpha+C_{\lambda}) C_{t^{*}}}
\bigg\{\Li{\lambda}_{1} C_{t^{*}} (C_{c}+\frac{C_{f}}{\alpha}) + \Big(\Li{c}_{1}+\Li{c}_{2}\Li{t^{*}}+C_{c}E_{1}\Big) \vee \Big( \Li{Q} \Li{w}_{*} \Big)\\
& & +2E_{3} +\frac{2\Li{Q}C_{\lambda}}{\alpha}\Li{w}_{1} + \Big\{E_{1}+2E_{2}+\Li{\lambda}_{1}C_{t^{*}}(1+C_{\lambda}/\alpha)\Big\} C_{w}
 \bigg\}, \\
\Li{\mathcal{L}w}_{2} & \leq & \mathrm{e}^{(\alpha+C_{\lambda}) C_{t^{*}}}
\Big\{ 3 C_{f} +\Li{c}_{2}+2 C_{c}(C_{\lambda}+\alpha) +\frac{C_{f}C_{\lambda}}{\alpha}
+C_{w} \big[ 4C_{\lambda}+ \frac{C_{\lambda}^{2}}{\alpha}+\alpha\big] \Big\}, \\
\Li{\mathcal{L}w}_{*} & \leq & \Li{\mathcal{L}w}_{1} + \Li{\mathcal{L}w}_{2} \Li{t}_{*}, \\
\Li{\mathcal{L}w} & \leq &  \big\{E_{1}+E_{2}\big\} C_{w} + \frac{\Li{Q}C_{\lambda}}{\alpha}\Li{w}_{1} + E_{3}
+ \Big(\Li{c}_{1}+\Li{c}_{2}\Li{t^{*}}+C_{c}E_{1}\Big) \vee \Big( \Li{Q} \Li{w}_{*} \Big).
\end{eqnarray*}
\end{proposition}

\noindent{\textbf{Proof:}} 
Let us denote $\mathcal{L}w$ by $g$.
We have for $(x,y)\in E^{2}$ and $t\in [0,\tstar{x}\wedge\tstar{y}]$ 
\begin{align*}
\Big| & g(\ph{x}{t})-g(\ph{y}{t}) \Big| \leq 
\mathrm{e}^{\alpha t+\Lambda(y,t)} \Big\{ \big| F(x,t)- F(y,t) \big| + \big| Iw(x,t)-Iw(y,t) \big| \Big\} \\
& + \mathrm{e}^{\alpha t+\Lambda(y,t)}
\Big\{   \big| \inf_{s\geq t} J(Mw,w) (x,s) - \inf_{s\geq t} J(Mw,w) (y,s) \big|  \vee  \big| Kw(x) - Kw(y) \big|  \Big\} \\
& +  \Big| \mathrm{e}^{\alpha t+\Lambda(x,t)}-\mathrm{e}^{\alpha t+\Lambda(y,t)} \Big|
\bigg| \Big\{ \inf_{s\geq t} J(Mw,w) (x,s) \wedge Kw(x) \Big\} -F(x,t)- Iw(x,t) \bigg|.
\end{align*}
It is easy to show that for $x\in E$, $t\in [0,\tstar{x}]$, and $w\in\mathbf{L}$ we have
$e^{\alpha t+\Lambda(x,t)}\leq\mathrm{e}^{(\alpha+C_{\lambda}) C_{t^{*}}}$, 
$\bigg| \Big\{ \inf_{s\geq t} J(Mw,w) (x,s) \wedge Kw(x) \Big\} -F(x,t)- Iw(x,t) \bigg| \leq \frac{1}{\alpha}(C_{f}+C_{\lambda}C_{w})+C_{c}+C_{w}$
and  for $(x,y)\in E^{2}$ and $t\in [0,\tstar{x}\wedge\tstar{y}]$ 
$ \Big| \mathrm{e}^{\alpha t+\Lambda(x,t)}-\mathrm{e}^{\alpha t+\Lambda(y,t)} \Big| \leq  \mathrm{e}^{(\alpha+C_{\lambda}) C_{t^{*}}} \Li{\lambda}_{1} C_{t^{*}} |x-y|$.
Consequently, by using Lemmas \ref{TechFI} and \ref{TechK} and Remark \ref{InfJx}, we get the first equation.

\bigskip

For $x \in E$ and $(s,t)\in [0,\tstar{x}]^{2}$
\begin{align*}
\Big| & g(\ph{x}{s})-g(\ph{x}{t}) \Big| \leq 
 \mathrm{e}^{\alpha t+\Lambda(x,t)} \big| Iw(x,s)-Iw(x,t) \big|  \\
& + \mathrm{e}^{\alpha t+\Lambda(x,t)}
\Big\{   \big| \inf_{u\geq s} J(Mw,w) (x,u) - \inf_{u\geq t} J(Mw,w) (x,u) \big|  + \big| F(x,s)- F(x,t) \big|  \Big\} \\
& + \Big| \mathrm{e}^{\alpha s+\Lambda(x,s)}-\mathrm{e}^{\alpha t+\Lambda(x,t)} \Big|
\bigg| \Big\{ \inf_{u\geq s} J(Mw,w) (x,u) \wedge Kw(x) \Big\} -F(x,s)- Iw(x,s) \bigg|.
\end{align*}
Note that for $x\in E$, $(s,t) \in [0,\tstar{x}]^{2}$
$ \Big| \mathrm{e}^{\alpha s+\Lambda(x,s)}-\mathrm{e}^{\alpha t+\Lambda(x,t)} \Big| 
\leq  \mathrm{e}^{(\alpha+C_{\lambda}) C_{t^{*}}} (C_{\lambda}+\alpha) |t-s|$.
Consequently, we have by using Lemmas \ref{TechFI} and \ref{TechInfJt}, we obtain the second equation.

\bigskip

The third inequality is straightforward and finally, for $(x,y)\in E^{2}$ we have
\begin{align*}
\Big| & g(x)-g(y) \Big| \leq \big| \inf_{s\geq 0} J(Mw,w) (x,s) - \inf_{s\geq 0} J(Mw,w) (y,s) \big|  \vee  \big| Kw(x) - Kw(y) \big|.
\end{align*}
By using Remark \ref{InfJx} and Lemma \ref{TechK}, we get the last equation.
\hfill $\Box$

\begin{corollary}\label{cor v lip}
For all $0\leq n\leq N$, the value functions $v_n$ are in $\mathbf{L}$.
\end{corollary}

\noindent{\textbf{Proof:}} 
As $v_N=g$ is in $\mathbf{L}$ by assumption, a recursive application of Proposition~\ref{TechLiL} yields the result.
\hfill $\Box$

\bibliographystyle{acm}
\bibliography{impulse}

\end{document}